\documentclass[10pt,onecolumn]{IEEEtran}

\usepackage{amsfonts, amsmath, amssymb}
\usepackage{array, tabularx}
\usepackage[utf8]{inputenc}
\usepackage{tikz}
\usepackage{cite}
\usepackage{pgfplots}
\pgfplotsset{compat=1.18}
\pgfplotsset{every tick label/.append style={font=\footnotesize}}
\pgfplotsset{every axis/.append style={font=\footnotesize}}

\definecolor{azulcito}{RGB}{0,80,110}
\definecolor{rojito}{RGB}{110,0,80}
\definecolor{verdecito}{RGB}{80,110,0}

\newtheorem{theorem}{Theorem}

\newtheorem{lemma}[theorem]{Lemma}
\newtheorem{proposition}[theorem]{Proposition}
\newtheorem{assumption}{Assumption}
\newtheorem{remark}{Remark}
\newtheorem{example}{Example}

\newcommand{\ip}[2]{\langle{#1},{#2}\rangle}

\begin{document}

\title{Dynamics and Optimization in Spatially Distributed Electrical Vehicle Charging}

\author{Fernando~Paganini, ~\IEEEmembership{Fellow,~IEEE,} 
and Andres~Ferragut~
\thanks{The authors are with the Mathematics Applied to Telecommunications and Energy Research Group, Universidad ORT Uruguay, Montevideo, Uruguay (e-mail: paganini@ort.edu.uy).} \thanks{Research supported by AFOSR-US under Grant FA9550-12-1-0398.} \thanks{Preprint submitted to IEEE Transactions on Control of Network Systems.}
}

\maketitle

\begin{abstract}
We consider a spatially distributed demand for electrical vehicle recharging, that must be covered by a fixed set of charging stations. Arriving EVs receive feedback on transport times to each station, and waiting times at congested ones, based on which they make a selfish selection. This selection determines total arrival rates in station queues, which are represented by a fluid state; departure rates are modeled under the assumption that clients have a given sojourn time in the system. The resulting differential equation system is analyzed with tools of optimization. We characterize the equilibrium as the solution to a specific convex program, which has connections to optimal transport problems, and also with road traffic theory. In particular a price of anarchy  appears with respect to a social planner's allocation. From a dynamical perspective, global convergence to equilibrium is established, with tools of Lagrange duality and Lyapunov theory. 
An extension of the model that makes customer demand elastic to observed delays is also presented, and analyzed with extensions of the optimization machinery. Simulations to illustrate the global behavior are presented, which also help validate the model beyond the fluid approximation. 
\end{abstract}


\section{Introduction}\label{sec.intro}

The application of optimization tools to operate spatially distributed
facilities has a very rich and extensive history. Perhaps the oldest such research is \emph{optimal transport}, that goes back to Monge in the 18th century (see, e.g. \cite{santambrogio2015}): the problem concerns a planner's decision on how to efficiently transport mass between demand and supply spatial distributions.

More recently, the operation of engineered \emph{networks} often resorts to optimization, e.g. for power dispatch in the electric grid \cite{momoh2017electric}, or routing in communication networks \cite{bertsekas2021data}. These decisions must \emph{dynamically} adapt to changing conditions and, as networks grow to enormous scale, \emph{distributed} decisions become imperative. Designs which meet these challenges may sometimes be found through a combination of differential equations and convex optimization machinery, as was the case for resource allocation in the Internet \cite{kelly1998rate, low2002internet}. 

A key question in distributed network operation is: to what degree may control be imposed on individual units or, on the contrary, are these agents making their own selfish decisions? A prominent instance of the latter case is routing in road traffic networks \cite{sheffi1985urban}: if properly informed, drivers naturally select paths of the least latency, which results in a congestion \emph{game}, whose equilibrium was classically analyzed by Wardrop \cite{wardrop1952road}. This solution, while not centrally planned, can nevertheless be characterized in terms of a suitable optimization problem, which has been helpful to understand the \emph{Price of Anarchy}, i.e. the gap between this equilibrium and the social welfare optimum \cite{roughgarden2005selfish}, and to propose means (e.g. tolls) to mitigate it. While much of this classical analysis concerns equilibrium, \emph{dynamic} studies of road traffic networks are also extensive, see e.g. \cite{como2021distributed} and references therein.

In this paper we consider a new application area, the operation of an  Electrical Vehicle (EV) charging infrastructure. In particular, we are interested in public facilities situated in parking lots, where EV chargers are made available for temporary use \cite{mozaffari2020joint}. This development has motivated an active area of research, within which we distinguish different problems: (i) the operation of a \emph{single} facility of this kind, in particular the scheduling of charging opportunities taking into account EV deadlines and installation limitations 
\cite{acnFramework,zeballos2019proportional};
(ii) integration of EV charging to the smart grid
\cite{mukherjee2014review,wang2016smart}; (iii) facility location problems, i.e. where to deploy EV charging \cite{kuby2005flow,PagEMF22}. 

Our focus here is on the \emph{operation} of a spatially \emph{distributed} infrastructure made up of several charging stations, to efficiently serve a distributed demand for EV recharging. At a high level, this appears to be an optimal transport problem: demand follows a certain spatial distribution, and supply is offered in another; the optimal allocation assigns demand to stations while minimizing the overall required travel. 

Charging demand does not, however, materialize in a batch; rather, we have a \emph{dynamic} situation in which requests for service arise asynchronously over time at different spatial locations. Selecting an adequate station for each request must consider both the transport cost and station congestion. Without the transport component, such \emph{load balancing} decisions have been studied extensively in computer networks, where routing is in charge of a central dispatcher (see e.g. \cite{der2022scalable}). Here we must incorporate the transport aspect and, importantly, compulsory routing is not assumed. 

Rather, drivers will select a station to obtain the fastest possible service, similarly to selfish routing for road traffic. Indeed, our results have some parallels with this literature, but also distinguishing features. The road traffic problem considers a network of links in which {all} traffic is affected by selfish routing, and delays are a static function of link flows\footnote{Some models also include vehicle \emph{densities} in links as variables \cite{como2015throughput,como2021distributed}, and sometimes PDE effects \cite{coclite2005traffic}; these are beyond our scope.}. Here we are analyzing an \emph{overlay} on the road system, in which a small portion of vehicles takes part; from this perspective, transport delays are exogenous, they depend on distance and background traffic but not on EV routing decisions. On the other hand, routing will affect congestion at charging stations; we assume congestion delay information is fed back, together with transport delays, to the selfish routing agents. 

Our contributions are as follows:
\begin{itemize}
    \item We present a differential equation model for the evolution of station queues, driven by spatially distributed demand for charge and selfish decisions on station assignment. A non-standard feature is the treatment of departures: rather than require that jobs remain until completion (i.e., complete charge), for public charging facilities it is most natural to assume that EVs depart after a given \emph{sojourn time}. This  assumption impacts both queue evolution models and the calculation of queueing delays. These waiting times, superimposed to the transport delays, determine the selfish routing flows and thus the full dynamics. 
    \item We characterize the dynamics by introducing a suitable \emph{convex optimization}, variant of the optimal transport problem: equilibrium points are proved to correspond to optima, and global convergence to equilibrium is established. Our results differ from other such characterizations in the selfish routing literature \cite{sheffi1985urban,roughgarden2005selfish}, due to our delay model as a function of the queue states. Proofs require extensive use of Lagrange duality, together with specialized refinements of Lyapunov-LaSalle theory. 
    \item We compare the resulting equilibrium with the socially optimal allocation, exhibiting the gaps between the two that lead to a price of anarchy. Again, while this coincides conceptually with classical selfish routing, there are differences in the specifics.  
    \item We extend the model to allow for \emph{elasticity} in the demand: as a function of the experienced delays, some drivers may choose not to participate in the recharge system. These modified dynamics are also related to optimization: introducing a utility function to model customer patience, we show that the equilibrium maximizes a certain surplus objective in the input rates. Proofs require the invocation of minimax theory for a suitably chosen convex-concave function. Global convergence is also established. 
    \item Simulation studies are carried out with a concrete instance of stations in the plane, to illustrate the resource allocation and its comparison with optimal transport. These experiments are based on   discrete, stochastic demands, demonstrating the approximate validity of our model beyond the fluid abstraction.
\end{itemize}
A preliminary version of some of our results appeared in the conference paper \cite{PagFAllerton23}. There we used a discontinuous, switching model for selfish routing, which made the dynamic study challenging. Here we employ a smooth approximation to switching, which allows for a complete mathematical treatment within the realm of ordinary differential equations. 

The rest of the paper is organized as follows. We first collect in Section \ref{ssec.prelim} some notation and background material. In Section \ref{sec.formu} we motivate the problem and develop our differential equation model. In Section \ref{sec.opt} we present the connection to convex optimization, and the resulting interpretations. Section \ref{sec.elastic}
covers the version with demand elasticity. Simulations are presented in Section \ref{sec.sims}, and conclusions given in Section \ref{sec.concl}. Some technical proofs are collected in the two Appendices. 

\subsection{Preliminaries and notation}\label{ssec.prelim}

We cover here some notation and background material from convex analysis. $\mathbb{R}^n$ is the standard $n$-dimensional space, and $\mathbb{R}_+^n=\{x\in \mathbb{R}^n: x_j \geq 0\ \forall j\}$ the non-negative orthant; 
$\Delta_n=\{\delta\in \mathbb{R}_+^n: \sum_{j=1}^n\delta_j=1\}$ is the unit simplex.

Matrix variables $X=(x_{ij})\in \mathbb{R}^{m\times n}$ are represented in uppercase, and we will use the notation $x^{i} = (x_{ij})_{j=1}^n\in \mathbb{R}^n$ to represent the $i$-th row vector of matrix $X$.

Both convex and concave functions $f:\mathbb{R}^n \to \mathbb{R}$ will appear; if they are differentiable, $\nabla f$ denotes the gradient. 
A basic non-smooth concave function $\varphi:\mathbb{R}^n\to \mathbb{R}$ is 
\begin{align}\label{eq.varphi inicial}
	\varphi(y) = \min_{j}(y_j) = \min_{\delta\in \Delta_n} \sum_{j=1}^n y_j\delta_j.
\end{align}

We will extensively use a smooth approximation to the minimum, a ``log-sum-exp" function with parameter $\epsilon > 0$:
\begin{align}\label{eq.logsumexp}
     \varphi_\epsilon(y):= -\epsilon \log\Bigg(\sum_j e^{-y_j/\epsilon}\Bigg).
\end{align}
This function may be called a ``soft-min", given the bounds: 
\begin{align}\label{eq.boundlogsumexp}
    \min(y_j)-\epsilon \log(n)\leq \varphi_\epsilon(y)\leq \min(y_j).
\end{align}
$\varphi_\epsilon(y)$ is concave, and its gradient $\nabla \varphi_\epsilon(y)=:\delta(y)$ is an element of the unit simplex $\Delta_n$, with components 
\begin{align}\label{eq.softargmin}
    \delta_{j}(y) = \frac{e^{-\frac{y_j}{\epsilon}}}{\sum_{k=1}^n e^{-\frac{y_k}{\epsilon}}};
\end{align}
these are largest for the minimizing coordinates of $y$. 

Introduce finally the \emph{negative entropy} function 
\begin{align}\label{eq.negentropy}
    \mathcal{H}(\delta) = \sum_{j=1}^n \delta_j \log(\delta_j), \quad \delta \in \Delta_n. 
\end{align}
$\mathcal{H}$ is strictly convex, non-positive and lower bounded over the unit simplex. It is connected to log-sum-exp
by convex (Fenchel) duality, as stated in the following (see e.g. 
\cite{boyd2004convex}):
\begin{lemma}\label{lem.fenchel}
    For $y\in \mathbb{R}_+^n$, 
    \[
    \varphi_\epsilon(y)=\min_{\delta \in \Delta_n} \Big[ \sum_{j=1}^n y_j\delta_j+\epsilon \mathcal{H}(\delta)\Big].
    \]
    Furthermore, the unique minimizing $\delta(y)$ is given by \eqref{eq.softargmin}. 
\end{lemma}

Note that as $\epsilon\to 0+$, the optimization  above approaches the one in \eqref{eq.varphi inicial}, consistently with the approximation \eqref{eq.boundlogsumexp}.

\section{Dynamic model}\label{sec.formu}

We consider a set of EV charging \emph{stations}, indexed by $j=1,\ldots,n$, occupying certain locations in the territory. We will not make explicit reference to this geometry, it will be encoded in the transport costs to be specified. Each station has capacity (number of charging spots, assumed all identical) $c_j$.

Vehicles demanding service are distributed in space. While conceivably there could be a continuum of locations from which a service request could arise, a natural simplification is to consider a discrete set of locations $i=1,\ldots m$; for instance, this could be the set of city corners, and we associate the demand location to the closest corner. Typically, $m\gg n$.

The relative positions of demand locations $i$ and stations $j$ are reflected in a matrix $\mathcal{K}=(\kappa_{ij})$ of transport costs. In this paper we assume cost has units of travel \emph{time}, and is exogenous, determined by distance and background traffic. This is an appropriate assumption if the vehicles participating in recharging are a minor portion of the overall traffic. 

The decision to be implemented is an \emph{assignment} of demand to supply, i.e. of charging requests to stations, that takes into account travel costs. As a starting point for our discussion, it is useful to consider first an \emph{optimal transport} problem\footnote{Here and henceforth, for brevity we will often omit the index ranges $i=1,\ldots m$, and $j=1,\ldots n$.}: 
\begin{gather}\label{eq.kant}
	\min \sum_{ij} \kappa_{ij}\xi_{ij}, \mbox{ subject to} \\
  \xi_{ij}\geq 0 \ \ \forall i,j; \quad \sum_j \xi_{ij}= {\eta}_i \ \ \forall i; \quad  \sum_i \xi_{ij}= {\sigma}_j \ \ \forall j. \nonumber
\end{gather}
In this discrete version of the Monge-Kantorovich \cite{santambrogio2015} problem\footnote{Also called the Hitchcock problem in the transportation literature \cite{sheffi1985urban}.}, $\eta_i$ and $\sigma_j$ would be given, respectively demand and supply quantities at each location, satisfying global balance ($\sum_j \sigma_j = \sum_i \eta_i$), and one seeks a \emph{transport plan} $\{\xi_{ij}\}$ of minimum cost. This is a \emph{static}, one shot allocation decision.

In our application, instead, demand is \emph{dynamic}: charging requests 
arrive over time, and are assigned upon arrival to a station, chosen with real-time information on the system. In addition to travel cost, current station congestion must be considered. Moreover, instead of a global planner we will have decentralized decisions consistent with selfish incentives. Our aim is to model this dynamic allocation process.

In our 
fluid model, the {rate} $r_i\geq 0$ of requests/sec arising at each location $i$ will be an exogenous quantity. Initially it will be fixed, later in Section \ref{sec.elastic} we will introduce elasticity in the demand process. The model will determine a matrix variable $X(t)=(x_{ij}(t))$, representing the rates of requests from location $i$ directed to station $j$, as a function of time. These will satisfy the balance conditions
\begin{align*}
x_{ij} \geq 0, \quad \sum_{j} x_{ij} = r_i \ \mbox{ for each }\ i=1,\ldots m,     
\end{align*}
of analogous form to the demand constraints in \eqref{eq.kant}, but here in units of \emph{rates}.  It is convenient to introduce also the variables $\delta_{ij}=x_{ij}/r_i$, denoting the \emph{fraction} of requests from location $i$ sent to station $j$. For each fixed $i$, $\delta^i =(\delta_{ij})_{j=1}^n$ is a vector in the $n$-dimensional unit simplex $\Delta_n$.

In contrast with optimal transport, here we do not specify the rates on the supply side, i.e. the EV/sec served at each station; this will be a consequence of routing decisions. If transportation costs were the only consideration, the natural assignment would be for EVs to choose the cheapest (closest in travel time) station. For a concrete visualization: if travel costs $\kappa_{ij}$ are proportional to Euclidean distance, this would mean breaking up demand spatially in accordance to the Voronoi cells (see e.g. \cite{boyd2004convex}) associated with station locations; Section \ref{sec.sims} has illustrative pictures.  

However, station congestion should also influence the allocation decision. At each station $j$, denote by $q_j$ the corresponding queue of EVs, whether or not they have a charging spot. In fluid terms, the queue dynamics are described by 
\begin{align}\label{eq.queue}
	\dot{q}_j = \sum_{i} x_{ij}-d_j, \quad j=1,\ldots,n;
\end{align}
here $d_j$ denotes the departure rate from the queue. Different models for departures could apply, in  
\cite{PagFAllerton23} we considered two options: either EVs depart when service is completed, as is standard in queueing theory, or they depart after a certain \emph{sojourn time}. In this paper we focus on the second option, more natural when dealing with public charging facilities.
Namely, customers have a time budget assigned to the recharge operation, and upon its expiration they will depart from the system, irrespective of the amount of service (recharge) they have obtained. 

In practice, drivers would have different sojourn times, which could be  modeled as random. For our fluid formulation, we will denote by $T$ the mean sojourn time across the population, and use the following expression 
 for the departure rate as a function of queue occupation:
\begin{align}\label{eq.dep sojourn}
	d_j(q_j) = \frac{q_j}{T}.
\end{align}
The above is a version of Little's law in queueing theory \cite{bertsekas2021data}, a basic conservation law between mean queues, rates and sojourn times. It is valid under very general conditions assuming balance between arrival and departure rates. Here we are applying it to departure only, allowing for a mismatch with arrivals during the transient regime. In Section \ref{sec.sims} we will test this transient model through simulations. 

In reference to \eqref{eq.queue}, note that since $d_j(0)=0$, and $x_{ij}\geq 0$, non-negativity of queues is automatically preserved.

\begin{remark}\label{rem.delay}
Eq. \eqref{eq.queue} does not distinguish between the time a vehicle arrives into the system and selects a station, and the time it effectively reaches it. Again our motivation is simplicity, otherwise we would have 
delay-differential equation; our simpler model is approximately valid provided transport delays are much smaller than the time-scale of sojourn times. We note that the \emph{equilibrium} characterization to be obtained below would not be affected by such delays, but convergence analysis would be far more difficult. 
\end{remark}

If the queue $q_j$ is below capacity, then all assigned vehicles have a corresponding slot, so there are no delays to receiving service, other than travel. On the contrary, if $q_j > c_j$ there will be an additional \emph{waiting time}, which we proceed to model, assuming a first-come first-serve queueing policy.  An arriving EV must wait for the time until the excess assignment $q_j - c_j$ is cleared by the departure process. Our model for waiting delay is thus\footnote{Henceforth, $[\cdot]^+ = \max\{\cdot,0\}$.}:
\begin{align}\label{eq.mu of q}
	\mu_j:= \frac{[q_j-c_j]^+}{d_j(q_j)} = T \left[1-\frac{c_j}{q_j}\right]^+,
\end{align}
where we have invoked \eqref{eq.dep sojourn}. The function $\mu_j(q_j)$ is depicted in Fig \ref{fig.mu} below.

\begin{figure}
	\centering
	\begin{tikzpicture}
		\begin{axis}[
			width=0.4\columnwidth,
			axis x line=center,
			axis y line=center,
			xtick={0,1},
			xticklabels={$0$,$c_j$},
			ytick={0},
			xlabel={$q_j$},
			xlabel style={below},
			ylabel style={left},
			xmin=-.1,
			xmax=2.5,
			ymin=-.02,
			ymax = .75]
		  \addplot [azulcito, very thick,mark=none,domain=0:1] {0};
		  \addplot [azulcito, very thick,mark=none,domain=1:2.5] {1-1/x};
  		  \node[azulcito] at (axis cs:1.2,0.6) {\footnotesize $\mu_j(q_j):= T \left[1-\frac{c_j}{q_j}\right]^+$};
		  \end{axis}
	\end{tikzpicture}\hspace{2em}%
	\begin{tikzpicture}
		\begin{axis}[
			width=0.4\columnwidth,
			axis x line=center,
			axis y line=center,
			xtick={0,1},
			xticklabels={$0$,$c_j$},
			ytick={0},
			xlabel={$q_j$},
			xlabel style={below},
			ylabel style={left},
			xmin=-.1,
			xmax=2.5,
			ymin=-.02]
		  \addplot [azulcito, very thick,mark=none,domain=0:1] {0};
		  \addplot [azulcito, very thick,mark=none,domain=1:2.5] {x-1-ln(x)};
		  \node[azulcito] at (axis cs:1.2,.5) {\footnotesize $\displaystyle \beta_j(q_j) = \frac{1}{T}\int_0^{q_j}\mu_j(\sigma) d\sigma$};
		\end{axis}
	\end{tikzpicture}

	\caption{Delay model \eqref{eq.mu of q} and penalty barrier function \eqref{eq.barrier function}.}\label{fig.mu}

\end{figure}
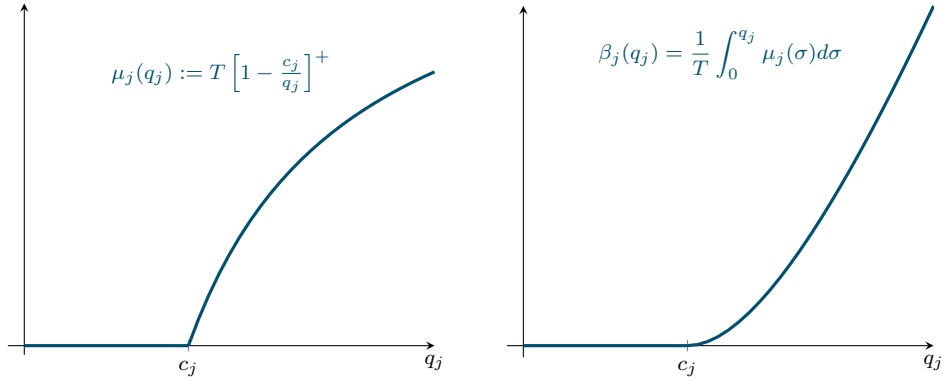

Consequently, an EV arriving at location $i$, if choosing station $j$, will be subject to a total delay to service of $\kappa_{ij}+\mu_j$. Note that  both terms have compatible units of \emph{time}.

The main assumption that completes our model is that the information $\kappa_{ij}+\mu_j$ is available to drivers, who then make a selfish decision.  Note, from a practical perspective, that time travel information is currently accessible through smartphone technology, and one could easily deploy an application through which stations broadcast their current waiting times. 

Let $y^i =\kappa^i + \mu$ be the vector of delays observed from location $i$ to all stations, i.e. $y^i_j = \kappa_{ij}+\mu_j, \ j=1,\ldots,n$. The natural selfish choice is $j\in \arg\min (y^i_j)$, i.e. minimizing delay to service.
We analyzed this model in \cite{PagFAllerton23}, which involves \emph{switching}, i.e., differential equations with a discontinuous field. While some of the analysis can be carried out in this setting, convergence proofs are technically challenging. 

In this paper we develop a smoother alternative, corresponding to the ``soft-min" approximation described in Section \ref{ssec.prelim}: routing fractions from location $i$ will follow the expression \eqref{eq.softargmin}, corresponding to  $y^i=\kappa^i + \mu$; smaller delays are favored, but in a less drastic fashion. Such ``logit" choice based on delay is a common model in the selfish routing literature (see \cite{sheffi1985urban}); it can be justified when there is noise in the delay information. Specifically, the routing fractions from location $i$ are:
\begin{align}\label{eq.softargmin-i}
    \delta_{ij}(\mathcal{K},\mu) = \frac{e^{-\frac{\kappa_{ij}+\mu_j}{\epsilon}}}{\sum_k e^{-\frac{\kappa_{ik}+\mu_k}{\epsilon}}}.
\end{align}
The preceding equation closes the feedback loop, leading to the overall dynamics: 
\begin{subequations}\label{eq.sojourn dynamics}
	\begin{align}
		\dot{q}_j =& \sum_{i=1}^m x_{ij} - \frac{q_j}{T}. \quad j=1,\ldots n.\label{eq.sojourn dynamics-state}\\
		\mu_j(q_j) =& T \left[1-\frac{c_j}{q_j}\right]^+,  \quad j=1,\ldots n. \label{eq.sojourn dynamics-mu} \\
  		x_{ij}=&r_i \delta_{ij},  \  i=1,\ldots m, j=1,\ldots n, \nonumber \\  &\mbox{ with } 
    \delta_{ij}(\mathcal{K},\mu)
    \mbox{ in \eqref{eq.softargmin-i}}. \label{eq.sojourn dynamics-x}
	\end{align}
\end{subequations}
The model has input parameters $T$, $c_j$, $r_i$, $\kappa_{ij}$. For the analysis to follow we assume constant transport costs, and thus omit henceforth the dependence on $\mathcal{K}$ in \eqref{eq.softargmin-i}.

We observe first that the above differential equation has a globally Lipschitz field. Indeed, the mapping $x_{ij}(\mu)$ is continuously differentiable in $\mathbb{R}^n$, and thus admits a global Lipschitz constant over the bounded set $\mu\in [0,T)^n$, which is the range of the function $\mu(q)$ in \eqref{eq.sojourn dynamics-mu}. The function $\mu_j(q_j)$ is not everywhere differentiable, but admits a global Lipschitz constant $T/c_j$. Thus $X(\mu(q))$ is globally Lipschitz; substitution into 
\eqref{eq.sojourn dynamics-state} gives a globally Lipschitz field in $q$. 

As a consequence, given an initial condition $q(0)$, solutions to  \eqref{eq.sojourn dynamics} exist, are unique, and defined for all time. In the following section, we will analyze their behavior invoking tools of convex optimization. 

\section{Optimization characterization} \label{sec.opt}

We begin by introducing a barrier function, which expresses a soft version of the capacity constraints. Let: 
\begin{equation}
	\beta_j(q_j) := \int_0^{q_j}\left[1-\frac{c_j}{\sigma}\right]^+ d\sigma = \begin{cases} 0 & q_j \leq c_j, \\
		q_j -c_j-c_j\log\left(\frac{q_j}{c_j}\right), & q_j > c_j. 
	\end{cases}\label{eq.barrier function}
\end{equation}
This is a convex, monotonically increasing function, also shown in Fig. \ref{fig.mu}.

We are now ready to introduce our convex optimization problem in the variables $X=(x_{ij})$ and $q=(q_j)$:
\begin{subequations}\label{eq.inelastic opt}
	\begin{gather}
		\min \sum_{i,j} \kappa_{ij}x_{ij} + \sum_j \beta_j(q_j) + \epsilon \sum_{i,j}x_{ij}\log\left(\frac{x_{ij}}{r_i}\right) \label{eq.costfc}
  \\
		\mbox{subject to:\quad}   x_{ij}\geq 0  \ \forall i,j; 
  \quad \sum_j x_{ij}= r_i \ \forall i; \label{eq.demand constr}\\ 
    \sum_i x_{ij} = \frac{{q}_j}{T}, \ \ \forall j. \label{eq.supply constr}
	\end{gather}
\end{subequations}
Let $C(X,q)$ be the cost function in \eqref{eq.costfc}. Its first term is analogous to the transportation cost in \eqref{eq.kant}; without the other terms, this amounts to the optimal transport of demands $r_i$ to supplies $q_j/T$, in units of (arrival and departure) rates. In contrast with \eqref{eq.kant}, however, the departure rates involve free decision variables $q_j$, penalized by the barrier cost \eqref{eq.barrier function}. 

An additional, relatively minor difference is the final perturbation term in the cost $C(X,q)$, which may be expressed as a negative entropy in the routing fractions:
\begin{align}
    \label{eq.negentropycost}
    \epsilon \sum_{i,j} r_i \delta_{ij} \log(\delta_{ij})
    = \epsilon\sum_i r_i \mathcal{H}(\delta^i);
\end{align}
this encourages diversity in the choice of station, slightly softening the selfish routing choice. 
\begin{remark}\label{rem.ri zero}
    The final term in \eqref{eq.costfc} is written for $r_i>0$; we could exclude a priori locations $i$ with $r_i=0$. In expression \eqref{eq.negentropycost} these terms drop out automatically: the penalty term goes to zero as demand disappears, a feature which will be relevant for the elastic case in the following section.  
    \end{remark}
\begin{proposition}\label{prop.unique opt}
The optimization problem \eqref{eq.inelastic opt} is feasible, and has a unique optimal point $(X^*,q^*)$.    
\end{proposition}
\begin{IEEEproof}
    For feasibility,  \eqref{eq.demand constr} can be satisfied by any split $x_{ij}$ of the demands $r_i>0$, and the free variable $q_j$ can accommodate \eqref{eq.supply constr}. For existence/uniqueness of the minimum, one can consider an equivalent problem in $X$ by replacing 
     $q_j$ from \eqref{eq.supply constr} into the cost term $\beta_j(\cdot)$. The resulting function of $X$ is \emph{strictly} convex (due to the negative entropy term), over a compact domain \eqref{eq.demand constr}: thus there is a unique minimum $X^*$, and consequently  $q^*$  due to \eqref{eq.supply constr}.      
\end{IEEEproof}

The main result of this section is that the optimization \eqref{eq.inelastic opt} characterizes the load balancing dynamics \eqref{eq.sojourn dynamics}: the optimal point is a globally attractive equilibrium of the dynamics. To prove this will require the use of Lagrange duality, as follows. 

\subsection{Lagrangian and Equilibrium Characterization}\label{ssec.lag-eq}
Introduce the Lagrangian of Problem  \eqref{eq.inelastic opt} with respect to  constraints \eqref{eq.supply constr}:

\begin{subequations}\label{eq.lagrangian}
\begin{align}
	L(X,q,\mu) =& \sum_{i,j} \kappa_{ij}x_{ij} +  \sum_j \beta_j(q_j)+  
 \epsilon \sum_{i,j}x_{ij}\log\left(\frac{x_{ij}}{r_i}\right) +\sum_{j} \mu_j \left[\sum_i x_{ij} - \frac{{q}_j}{T}\right]	\label{eq.lag}\\
	= &\underbrace{\sum_{i,j} x_{ij}\Big [\kappa_{ij}+\mu_j + \epsilon\log\left(\frac{x_{ij}}{r_i}\right)\Big]}_{\mbox{\normalsize $L_1(X,\mu)$}} \label{eq.lag1}\\
& + \underbrace{\sum_j \left[\beta_j(q_j)-\mu_j\frac{{q}_j}{T}\right]}_{\mbox{\normalsize $L_2(q,\mu)$}}. \label{eq.lag2}
\end{align}
\end{subequations}

Suggestively, we have denoted the multipliers by $\mu_j$; the optimum of the convex program \eqref{eq.inelastic opt} will correspond to a saddle point of this Lagrangian (minimum in $(X,q)$, maximum in $\mu$). We state the following result. 

\begin{theorem}\label{teo.sojourn-eq} 
	The following are equivalent:
\begin{itemize}
\item[(i)] $(X^*,q^*,\mu^*)$ is the saddle point of the Lagrangian $L$ in \eqref{eq.lagrangian}.
\item[(ii)] $(X^*,q^*,\mu^*)$ is an equilibrium point of 
\eqref{eq.sojourn dynamics}, under constant $r_i$.   
\end{itemize}
In particular, the dynamics have a unique equilibrium point.
\end{theorem}
\medskip
\begin{IEEEproof}
First observe, focusing on \eqref{eq.lag}, that  
for the maximum over an unconstrained $\mu$ to be finite requires primal feasibility of 
\eqref{eq.supply constr}, i.e., equilibrium of \eqref{eq.sojourn dynamics-state}:
\begin{align}\label{eq.saddle mu}
  \mu^* \in \arg\max_\mu L(X^*,q^*,\mu) \ \Longleftrightarrow \ \sum_i x^*_{ij} = \frac{q^*_j}{T},\; j=1,\ldots n.  
\end{align}
We now look at the minimization over $(X,q)$, treating both terms 
in \eqref{eq.lag1} and \eqref{eq.lag2} separately. 

For the minimization of $L_1(X,\mu)$ we must consider the remaining constraints \eqref{eq.demand constr} on $X$, which decouple across $i$. Invoking the routing fractions $\delta_{ij}=x_{ij}/r_i$ we write
\begin{align}\label{eq.lag1 gamma}
    L_1(X,\mu)= \sum_i r_i \left[ \sum_j \delta_{ij}(\kappa_{ij}+\mu_j) + \epsilon \mathcal{H}(\delta^i)\right].
\end{align}
To minimize each term in square brackets over the unit simplex, we apply 
Lemma \ref{lem.fenchel} for the vector $y^i=\kappa^i + \mu$. The result is the log-sum-exp expression\footnote{This new notation emphasizes the dependence on the variable $\mu$; the super-index $i$ in $\varphi^i_\epsilon(\mu)$ indicates the displacement of $\mu$ by the fixed vector of transport costs $\kappa^i$ from location $i$.}
\begin{align}\label{eq.logsumexp-i}
\varphi^{i}_\epsilon(\mu):=\varphi_\epsilon(\kappa^i+\mu)= -\epsilon \log\left(\sum_j e^{-(\kappa_{ij}+\mu_j)/\epsilon}\right),
\end{align}
achieved at $\delta^i(\mu)=(\delta_{ij}(\mu))_{j=1}^n\in \Delta_n$, which follows precisely the expression \eqref{eq.softargmin-i}. We thus conclude that 
\begin{equation}\label{eq.saddle X}
  X^* \in \arg\max_X L_1(X,\mu^*) \ \Longleftrightarrow \ x^*_{ij} = r_i \delta_{ij}(\mu^*), \mbox{ with }\delta_{ij}(\mu) \mbox{ in \eqref{eq.softargmin-i}}.  \nonumber 
\end{equation}

Finally, we consider the minimization of $L_2(\mu,q)$ which is unconstrained in $q$, and decoupled over $j$; we minimize each term separately.
Let 
\begin{align}
    D_{2j}(\mu_j)= \inf_{q_j} [\beta_j(q_j)-\mu_j{q}_j/T].
\end{align}
We have the following cases:
\begin{itemize} 
    \item If $\mu_j < 0$, or $\mu_j \geq T$, $ D_{2j}(\mu_j)=-\infty$.
    \item If $\mu_j = 0$, $ D_{2j}(\mu_j)=0$, achieved at $q_j\in (-\infty,c_j]$.
    \item If $0<\mu_j<T$, the unique minimizing $q_j$ is obtained by:
    \begin{align}\label{eq.phi prime inverse}    
      \beta_j'({q}_j)= \left[1-\frac{c_j}{q_j}\right] = \frac{\mu_j}{T} \ 
     \Longleftrightarrow \     q_j = \frac{T c_j}{T-\mu_j}.
\end{align}
\end{itemize}
We can encompass the two cases with finite minimum by the relationship
$\mu_j= T \left[1-\frac{c_j}{q_j}\right]^+$ between $\mu_j\in [0,T)$ and the minimizing $q_j$. From here we conclude that
\begin{align}\label{eq.saddle q}
  q^* \in \arg\max_q L_2(q,\mu^*) \ \Longleftrightarrow \ &\mu_j^*= T \left[1-\frac{c_j}{q_j^*}\right]^+.   
\end{align}
The left-hand sides of equations \eqref{eq.saddle mu},  \eqref{eq.saddle X}, and \eqref{eq.saddle q} are the saddle point conditions (i). The corresponding right-hand sides are the equilibrium conditions (ii) for the dynamics \eqref{eq.sojourn dynamics}.  

For the last statement, since the primal variables $(X^*,q^*)$ at a saddle point are optima, we can invoke Proposition \eqref{prop.unique opt} to show they are unique. Note that $\mu(q^*)$ from \eqref{eq.sojourn dynamics-mu} must be unique as well. 
\end{IEEEproof}

\subsection{Dual function} \label{ssec.dual}
As additional conclusion of the preceding analysis, let us make explicit the dual function 
\[
D(\mu)=\inf_{X,q}L(X,q,\mu) = \underbrace{\inf_X L_1(X,\mu)}_{D_1(\mu)} + \underbrace{\inf_q L_2(q,\mu)}_{D_2(\mu)}. 
\]
Referring back to \eqref{eq.lag1 gamma} and \eqref{eq.logsumexp-i}, we have 
\begin{align*}
	D_1(\mu) & = \min_{X \in \eqref{eq.demand constr}} L_1(X,\mu) = \sum_i r_i\varphi^{i}_\epsilon(\mu);
\end{align*}
we further note that $D_1(\mu)$ is a differentiable function 
of $\mu \in \mathbb{R}^n$, and its gradient takes the form:
\begin{align}\label{eq.d1grad}
	\nabla{D_1}(\mu) =   \sum_i r_i \nabla \varphi^{i}_\epsilon(\mu)
 = \sum_i r_i \delta^i(\mu) = \sum_i x^i(\mu),
\end{align}
where the last expression is based on \eqref{eq.sojourn dynamics-x}.

Secondly, in reference to \eqref{eq.phi prime inverse},  substitution of the minimizing $q_j$ into $[\beta_j(q_j)-\mu_j{q}_j/T]$ gives a minimum of
\begin{align}\label{eq.d2j}
    D_{2j}(\mu_j)=c_j\log(1-\mu_j/T),
\end{align}
and this formula also covers the case $\mu_j=0$. Therefore:
\begin{align}\label{eq.dual2}
	D_{2}(\mu) = \sum_{j} c_j \log\left(1-\frac{\mu_j}{T}\right),\quad \mbox{ for } 0 \leq \mu_j < T.  
\end{align}
This function is differentiable in the interior of the domain, but the boundary $\mu_j=0$ requires some special care, as we will see in the convergence analysis below. 

The overall dual function, with domain $\mu\in [0,T)^n$, is:
\begin{align}\label{eq.dual}
D(\mu) & = \sum_i r_i\varphi^{i}_\epsilon(\mu)+\sum_j  c_j \log\left(1-\frac{\mu_j}{T}\right).
\end{align}

\begin{proposition}\label{prop.dual}
    $D(\mu)$ is strictly concave and has a finite maximum $D^*$ over $\mu\in [0,T)^n$, achieved at a unique $\mu^*\in [0,T)^n$. 
\end{proposition}
\begin{IEEEproof}
    $D_1(\mu)$ is concave (not strictly), and $D_2(\mu)$ is strictly concave in $[0,T)^n$ (this follows directly from each component $D_{2j}$). Thus, $D(\mu)$ is strictly concave. Since  $D_{2j}(\mu_j) \to -\infty$ as $\mu_j \uparrow T$, there is a global maximum $D^* = D(\mu^*)$ with $\mu^*$ strictly within  $[0,T)^n$, unique due to strict concavity. This is also consistent with the uniqueness of $\mu^*$ in the saddle point, shown in Theorem \ref{teo.sojourn-eq}.
\end{IEEEproof}

\subsection{Interpretation and Price of Anarchy} \label{ssec.eq-interpr}

Equilibrium points for our model of dynamic station assignment have been shown to be solutions of a certain modified optimal transport problem. We now provide some interpretations of the result, and connections to the selfish routing literature \cite{sheffi1985urban,roughgarden2005selfish}. For simplicity, we will ignore in the discussion the entropy regularization term in the cost and focus on 
\begin{align}\label{eq.c0}
 C_0(X,q) := \sum_{i,j} \kappa_{ij}x_{ij} + \sum_j \beta_j(q_j);  
\end{align}
our equilibrium $(X^*,q^*)$ optimizes (approximately as $\epsilon \to 0$) this convex function, subject to constraints \eqref{eq.demand constr}-\eqref{eq.supply constr}. 

The first term above has a natural interpretation. Recall that $\kappa_{ij}$ represents the travel times, and $x_{ij}$ the rates, between arrival location $i$ and station $j$. Therefore, $\sum_{i,j} \kappa_{ij}x_{ij}$ will be the total number of EVs currently in travel towards a charging station, a natural transport cost to be minimized. 

The second term, based on $\beta_j(\cdot)$ in \eqref{eq.barrier function}, does not have such a transparent form. From a \emph{social welfare} perspective, the \emph{congestion cost} to add would be the total number of EVs waiting at stations without a charging spot: $\sum_j  \left[q_j-{c_j}\right]^+$.
The natural social welfare cost would thus be:
\begin{align}\label{eq.cs}
 C_s(X,q)  &= \sum_{i,j} \kappa_{ij}x_{ij} + \sum_j  \left[q_j-{c_j}\right]^+.
\end{align}
We remark the following:
\begin{itemize}
\item $C_s(X,q)$ is also convex; in fact it is piecewise linear\footnote{We could obtain \emph{strict} convexity over the domain, and thus a unique optimal point for
$C_s(X,q)$, by adding an entropy term as in \eqref{eq.costfc}.}.  Using slack variables $z_j$, the social planner's solution to the station assignment problem can be solved by means of the linear program: 
\begin{align*}
    \min& \ \kappa_{ij}x_{ij} + \sum_j  z_j, \\
    \mbox{ subjet to }& \mbox{\eqref{eq.demand constr}, \eqref{eq.supply constr}}, z_j \geq 0, \ z_j \geq q_j - c_j.
\end{align*}
\item $C_s(X,q) \geq C_0(X,q)$. Indeed, by \eqref{eq.barrier function} we have 
\begin{equation*}
C_s(X,q) - C_0(X,q) = \sum_j \left[q_j-{c_j}\right]^+ - \beta_j(q_j) =  \sum_j c_j\log\left(\frac{q_j}{c_j}\right) \mathbf{1}_{q_j > c_j}.
\end{equation*}
\end{itemize}

Note that the right-hand side above is zero in the absence of congestion 
($q_j \leq c_j \ \forall j$). A consequence is that if the equilibrium $(X^*,q^*)$ of our dynamics (which minimizes $C_0$, again ignoring the entropy term) involves no congestion, then it must also minimize 
$C_s(X,q)$: selfish routing achieves global welfare in this case. 

However, if station congestion appears, there is a difference between both costs, the equilibrium $(X^*,q^*)$ will no longer be socially optimal. This \emph{Price of Anarchy} appears for similar reasons as in road traffic models \cite{roughgarden2005selfish}. To see this, rewrite the social congestion cost at station $j$ (using \eqref{eq.mu of q})  as:
\[
 \left[q_j-{c_j}\right]^+ =\sum_j  \left[1-\frac{c_j}{q_j}\right]^+ q_j
=  \frac{1}{T}\sum_j \mu_j q_j,
 \]
and compare it to the corresponding barrier term in \eqref{eq.c0}:
\[
\beta_j(q_j)=\frac{1}{T}\int_0^{q_j} \mu_j(\sigma)d\sigma.
\]
An analogous variation appears in the road traffic literature, where latencies are taken to be static functions of link \emph{flows}. The product of latency and flow is the natural welfare cost, whereas the integral of the latency function is the cost that characterizes the Wardrop equilibrium in terms of optimization. 

The difference here is that our congestion delays appear at stations rather than links, depending on \emph{queues}, not flows; queues are state variables, the output of an integrator driven by flows. For this reason our results do not follow from the classical work, they require specific mathematics, in particular the extensive use of Lagrange duality, for both the equilibrium analysis and the convergence results below.

\begin{remark}
This distinction has parallels with the Internet congestion control literature \cite{kelly1998rate,low2002internet}. In that case, so-called ``primal" models take queues to be static functions of link flows, whereas ``dual" models favor the fluid integrator. The terminology comes from its association with convex optimization: primal models are analyzed in the original flow variables, dual models in the space of multipliers. 
\end{remark}

\begin{example}\label{example:poa}
We illustrate the Price of Anarchy through a toy example, solved numerically. There are two stations with capacities $c_1=20$, $c_2=40$, and a single source location. Transport times are  $\kappa_1=1, \kappa_2=10$, sojourn time is $T=60$. We increment the demand $r$ upward from zero, and compare the equilibrium with the socially optimal solution. 

In this simple case, $x_i=q_i/T$ so we can write the social welfare cost as a function of rates:
\[
C_s(x_1,x_2) = \kappa_1 x_1 + \kappa_2 x_2 +  \left[Tx_1-{c_1}\right]^+ +  \left[Tx_2-{c_2}\right]^+.
\]

In Fig. \ref{fig.poarates} (left) we plot the optimal breakup of rates as a function of $r$;
initially, $x_1=r$ so all traffic is sent to the closest station. 
Upon reaching $r=1/3$  this station fills ($rT=c_1$); and after that it becomes optimal to set $x_2=r-1/3$, i.e. send excess traffic to station 2. This is clear from the expression above: since $\kappa_2<T$, it is cheaper to pay the transport cost to station 2, than the congestion cost at station 1. 
After $r=1$, both stations are full, so congestion cost is inevitable and also indifferent to station choice. Therefore, travel costs dictate that further increases in traffic must be sent again to station 1. 

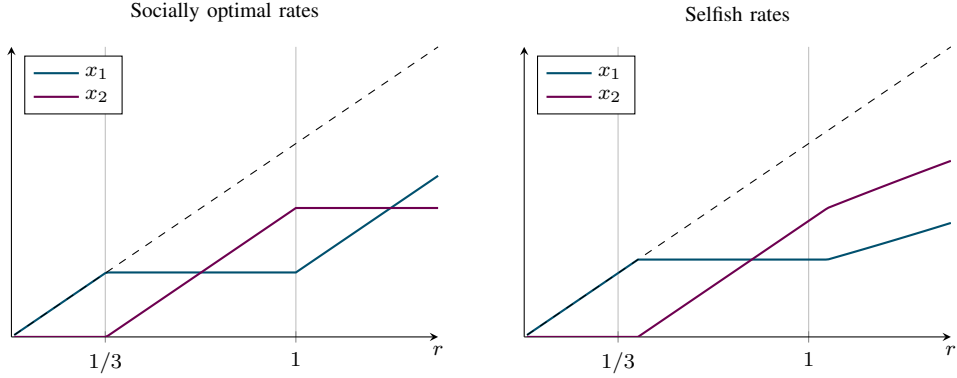
\begin{figure}
    \centering
    \begin{tikzpicture}
        \begin{axis}[
            width=0.4\textwidth,
            height=0.3\textwidth,
            title={Socially optimal rates},
            legend style={font=\footnotesize},
            axis x line=center,
            axis y line=center,
            xtick={0,0.33,1},
            xticklabels={$0$,$1/3$, $1$},
            ytick={0},
            xlabel={$r$},
            xlabel style={below},
            ylabel style={left},
            xmin=0,
            xmax=1.5,
            ymin=0,
            ymax = 1.5,
            legend pos = north west,
            grid=both]
        \addplot[color=azulcito, thick, solid]
            table[row sep={\\}]
            {
                \\
                0.01  0.010000000333317586  \\
                0.02  0.020000000333317593  \\
                0.03  0.03000000033331758  \\
                0.04  0.0400000003333176  \\
                0.05  0.05000000033331759  \\
                0.06  0.06000000033331759  \\
                0.07  0.0700000003333176  \\
                0.08  0.0800000003333176  \\
                0.09  0.09000000033331758  \\
                0.1  0.1000000003333176  \\
                0.11  0.1100000003333176  \\
                0.12  0.12000000033331759  \\
                0.13  0.1300000003333176  \\
                0.14  0.1400000003333176  \\
                0.15  0.1500000003333176  \\
                0.16  0.1600000003333176  \\
                0.17  0.1700000003333176  \\
                0.18  0.18000000033331756  \\
                0.19  0.1900000003333176  \\
                0.2  0.2000000003333176  \\
                0.21  0.2100000003333176  \\
                0.22  0.22000000033331762  \\
                0.23  0.2300000003333176  \\
                0.24  0.2400000003333176  \\
                0.25  0.2500000003333176  \\
                0.26  0.2600000003333176  \\
                0.27  0.2700000003333176  \\
                0.28  0.28000000033331757  \\
                0.29  0.2900000003333176  \\
                0.3  0.3000000003333176  \\
                0.31  0.31000000033331765  \\
                0.32  0.3200000003333176  \\
                0.33  0.3300000003333176  \\
                0.34  0.33333333649998853  \\
                0.35  0.33333333649998853  \\
                0.36  0.3333333364999885  \\
                0.37  0.33333333649998853  \\
                0.38  0.3333333364999885  \\
                0.39  0.33333333649998853  \\
                0.4  0.33333333649998853  \\
                0.41  0.33333333649998853  \\
                0.42  0.3333333364999885  \\
                0.43  0.3333333364999885  \\
                0.44  0.33333333649998853  \\
                0.45  0.33333333649998853  \\
                0.46  0.33333333649998853  \\
                0.47  0.33333333649998853  \\
                0.48  0.3333333364999885  \\
                0.49  0.33333333649998853  \\
                0.5  0.33333333649998853  \\
                0.51  0.33333333649998853  \\
                0.52  0.3333333364999885  \\
                0.53  0.3333333364999885  \\
                0.54  0.33333333649998853  \\
                0.55  0.3333333364999885  \\
                0.56  0.3333333364999885  \\
                0.57  0.3333333364999885  \\
                0.58  0.3333333364999885  \\
                0.59  0.33333333649998853  \\
                0.6  0.33333333649998853  \\
                0.61  0.33333333649998853  \\
                0.62  0.33333333649998853  \\
                0.63  0.3333333364999885  \\
                0.64  0.3333333364999885  \\
                0.65  0.3333333364999885  \\
                0.66  0.3333333364999885  \\
                0.67  0.3333333364999885  \\
                0.68  0.3333333364999885  \\
                0.69  0.3333333364999885  \\
                0.7  0.33333333649998853  \\
                0.71  0.3333333364999885  \\
                0.72  0.33333333649998853  \\
                0.73  0.3333333364999885  \\
                0.74  0.33333333649998853  \\
                0.75  0.33333333649998853  \\
                0.76  0.3333333364999885  \\
                0.77  0.33333333649998853  \\
                0.78  0.3333333364999885  \\
                0.79  0.3333333364999885  \\
                0.8  0.3333333364999885  \\
                0.81  0.33333333649998853  \\
                0.82  0.33333333649998853  \\
                0.83  0.3333333364999885  \\
                0.84  0.33333333649998853  \\
                0.85  0.3333333364999885  \\
                0.86  0.3333333364999885  \\
                0.87  0.3333333364999885  \\
                0.88  0.3333333364999885  \\
                0.89  0.3333333364999885  \\
                0.9  0.33333333649998853  \\
                0.91  0.3333333364999885  \\
                0.92  0.3333333364999885  \\
                0.93  0.3333333364999885  \\
                0.94  0.3333333364999885  \\
                0.95  0.3333333364999885  \\
                0.96  0.3333333364999885  \\
                0.97  0.33333333649998853  \\
                0.98  0.3333333364999885  \\
                0.99  0.33333333649998853  \\
                1.0  0.33333333649998853  \\
                1.01  0.3433333268333219  \\
                1.02  0.353333326833322  \\
                1.03  0.36333332683332187  \\
                1.04  0.373333326833322  \\
                1.05  0.38333332683332194  \\
                1.06  0.39333332683332195  \\
                1.07  0.40333332683332196  \\
                1.08  0.41333332683332186  \\
                1.09  0.42333332683332203  \\
                1.1  0.433333326833322  \\
                1.11  0.4433333268333221  \\
                1.12  0.4533333268333219  \\
                1.13  0.4633333268333219  \\
                1.14  0.47333332683332185  \\
                1.15  0.48333332683332186  \\
                1.16  0.4933333268333218  \\
                1.17  0.5033333268333218  \\
                1.18  0.5133333268333219  \\
                1.19  0.5233333268333219  \\
                1.2  0.5333333268333219  \\
                1.21  0.5433333268333218  \\
                1.22  0.5533333268333217  \\
                1.23  0.5633333268333219  \\
                1.24  0.5733333268333218  \\
                1.25  0.5833333268333218  \\
                1.26  0.593333326833322  \\
                1.27  0.6033333268333219  \\
                1.28  0.6133333268333219  \\
                1.29  0.6233333268333219  \\
                1.3  0.6333333268333219  \\
                1.31  0.6433333268333219  \\
                1.32  0.653333326833322  \\
                1.33  0.663333326833322  \\
                1.34  0.673333326833322  \\
                1.35  0.683333326833322  \\
                1.36  0.693333326833322  \\
                1.37  0.703333326833322  \\
                1.38  0.7133333268333217  \\
                1.39  0.7233333268333217  \\
                1.4  0.7333333268333216  \\
                1.41  0.7433333268333217  \\
                1.42  0.7533333268333218  \\
                1.43  0.7633333268333218  \\
                1.44  0.7733333268333217  \\
                1.45  0.7833333268333218  \\
                1.46  0.7933333268333218  \\
                1.47  0.8033333268333219  \\
                1.48  0.8133333268333218  \\
                1.49  0.8233333268333219  \\
                1.5  0.8333333268333221  \\
            }
            ;
        \addlegendentry {$x_1$}
        \addplot[color=rojito, thick, solid]
            table[row sep={\\}]
            {
                \\
                0.01  -3.3331758534761013e-10  \\
                0.02  -3.3331759057672884e-10  \\
                0.03  -3.3331758047089874e-10  \\
                0.04  -3.3331759490278573e-10  \\
                0.05  -3.333175883555897e-10  \\
                0.06  -3.333175917126786e-10  \\
                0.07  -3.3331759527464064e-10  \\
                0.08  -3.3331759933273393e-10  \\
                0.09  -3.3331757990702073e-10  \\
                0.1  -3.333175864215268e-10  \\
                0.11  -3.3331759452224005e-10  \\
                0.12  -3.3331759364629285e-10  \\
                0.13  -3.3331759323821436e-10  \\
                0.14  -3.3331759416314456e-10  \\
                0.15  -3.333175939698777e-10  \\
                0.16  -3.333175791885869e-10  \\
                0.17  -3.3331758375888083e-10  \\
                0.18  -3.3331758792864923e-10  \\
                0.19  -3.3331759155644326e-10  \\
                0.2  -3.3331759466376704e-10  \\
                0.21  -3.333175973144848e-10  \\
                0.22  -3.3331759683498137e-10  \\
                0.23  -3.3331759791697367e-10  \\
                0.24  -3.333175988371871e-10  \\
                0.25  -3.333175996081459e-10  \\
                0.26  -3.333176002424641e-10  \\
                0.27  -3.3331760000098104e-10  \\
                0.28  -3.333175723505787e-10  \\
                0.29  -3.333175997883743e-10  \\
                0.3  -3.3331760765382125e-10  \\
                0.31  -3.3331761688316214e-10  \\
                0.32  -3.333176128871888e-10  \\
                0.33  -3.333176058919866e-10  \\
                0.34  0.006666663500011528  \\
                0.35  0.016666663500011458  \\
                0.36  0.02666666350001149  \\
                0.37  0.0366666635000115  \\
                0.38  0.04666666350001152  \\
                0.39  0.05666666350001148  \\
                0.4  0.06666666350001146  \\
                0.41  0.07666666350001146  \\
                0.42  0.08666666350001151  \\
                0.43  0.0966666635000115  \\
                0.44  0.10666666350001146  \\
                0.45  0.1166666635000115  \\
                0.46  0.1266666635000115  \\
                0.47  0.13666666350001147  \\
                0.48  0.1466666635000115  \\
                0.49  0.15666666350001146  \\
                0.5  0.16666666350001147  \\
                0.51  0.17666666350001145  \\
                0.52  0.18666666350001146  \\
                0.53  0.1966666635000115  \\
                0.54  0.20666666350001156  \\
                0.55  0.21666666350001157  \\
                0.56  0.2266666635000116  \\
                0.57  0.23666666350001145  \\
                0.58  0.24666666350001146  \\
                0.59  0.2566666635000114  \\
                0.6  0.2666666635000114  \\
                0.61  0.27666666350001146  \\
                0.62  0.28666666350001146  \\
                0.63  0.29666666350001153  \\
                0.64  0.3066666635000116  \\
                0.65  0.31666666350001155  \\
                0.66  0.32666666350001156  \\
                0.67  0.33666666350001156  \\
                0.68  0.3466666635000116  \\
                0.69  0.35666666350001147  \\
                0.7  0.3666666635000114  \\
                0.71  0.3766666635000115  \\
                0.72  0.3866666635000114  \\
                0.73  0.39666666350001145  \\
                0.74  0.40666666350001146  \\
                0.75  0.41666666350001147  \\
                0.76  0.42666666350001153  \\
                0.77  0.4366666635000115  \\
                0.78  0.4466666635000115  \\
                0.79  0.45666666350001156  \\
                0.8  0.46666666350001157  \\
                0.81  0.4766666635000115  \\
                0.82  0.4866666635000115  \\
                0.83  0.49666666350001143  \\
                0.84  0.5066666635000115  \\
                0.85  0.5166666635000114  \\
                0.86  0.5266666635000115  \\
                0.87  0.5366666635000115  \\
                0.88  0.5466666635000116  \\
                0.89  0.5566666635000115  \\
                0.9  0.5666666635000115  \\
                0.91  0.5766666635000116  \\
                0.92  0.5866666635000115  \\
                0.93  0.5966666635000116  \\
                0.94  0.6066666635000114  \\
                0.95  0.6166666635000115  \\
                0.96  0.6266666635000114  \\
                0.97  0.6366666635000114  \\
                0.98  0.6466666635000115  \\
                0.99  0.6566666635000115  \\
                1.0  0.6666666635000115  \\
                1.01  0.666666673166678  \\
                1.02  0.6666666731666782  \\
                1.03  0.666666673166678  \\
                1.04  0.6666666731666782  \\
                1.05  0.666666673166678  \\
                1.06  0.6666666731666782  \\
                1.07  0.666666673166678  \\
                1.08  0.6666666731666782  \\
                1.09  0.6666666731666782  \\
                1.1  0.6666666731666782  \\
                1.11  0.666666673166678  \\
                1.12  0.6666666731666782  \\
                1.13  0.666666673166678  \\
                1.14  0.6666666731666782  \\
                1.15  0.666666673166678  \\
                1.16  0.666666673166678  \\
                1.17  0.6666666731666782  \\
                1.18  0.666666673166678  \\
                1.19  0.666666673166678  \\
                1.2  0.666666673166678  \\
                1.21  0.6666666731666782  \\
                1.22  0.6666666731666782  \\
                1.23  0.666666673166678  \\
                1.24  0.666666673166678  \\
                1.25  0.666666673166678  \\
                1.26  0.666666673166678  \\
                1.27  0.6666666731666782  \\
                1.28  0.6666666731666782  \\
                1.29  0.6666666731666782  \\
                1.3  0.6666666731666782  \\
                1.31  0.666666673166678  \\
                1.32  0.666666673166678  \\
                1.33  0.666666673166678  \\
                1.34  0.666666673166678  \\
                1.35  0.666666673166678  \\
                1.36  0.666666673166678  \\
                1.37  0.666666673166678  \\
                1.38  0.6666666731666782  \\
                1.39  0.6666666731666782  \\
                1.4  0.6666666731666782  \\
                1.41  0.6666666731666782  \\
                1.42  0.6666666731666782  \\
                1.43  0.6666666731666782  \\
                1.44  0.6666666731666782  \\
                1.45  0.6666666731666782  \\
                1.46  0.666666673166678  \\
                1.47  0.666666673166678  \\
                1.48  0.666666673166678  \\
                1.49  0.666666673166678  \\
                1.5  0.666666673166678  \\
            }
            ;
        \addlegendentry {$x_2$}
        
        \addplot[color=black, dashed, forget plot]
            table[row sep={\\}]
            {
                \\
                0.01  0.01  \\
                1.5  1.5  \\
            }
            ;
    \end{axis}
    \end{tikzpicture}
    \hspace{2em}
    \begin{tikzpicture}
        \begin{axis}[
            width=0.4\textwidth,
            height=0.3\textwidth,
            title={Selfish rates},
            axis x line=center,
            axis y line=center,
            xtick={0,0.33,1},
            xticklabels={$0$,$1/3$, $1$},
            ytick={0},
            xlabel={$r$},
            xlabel style={below},
            ylabel style={left},
            xmin=0,
            xmax=1.5,
            ymin=0,
            ymax = 1.5,
            legend pos = north west,
            grid = both]
    
        ]
        \addplot[color=azulcito, thick, solid]
            table[row sep={\\}]
            {
                \\
                0.01  0.010000000333315152  \\
                0.02  0.020000000333308246  \\
                0.03  0.030000000333308244  \\
                0.04  0.04000000033330824  \\
                0.05  0.050000000333308245  \\
                0.06  0.060000000333315144  \\
                0.07  0.07000000033331516  \\
                0.08  0.08000000033331515  \\
                0.09  0.09000000033330824  \\
                0.1  0.10000000033330825  \\
                0.11  0.11000000033330824  \\
                0.12  0.12000000033330824  \\
                0.13  0.13000000033330825  \\
                0.14  0.14000000033330826  \\
                0.15  0.15000000033330824  \\
                0.16  0.16000000033330825  \\
                0.17  0.17000000033330825  \\
                0.18  0.18000000033330824  \\
                0.19  0.19000000033330824  \\
                0.2  0.20000000033330825  \\
                0.21  0.21000000033330823  \\
                0.22  0.22000000033330824  \\
                0.23  0.23000000033330825  \\
                0.24  0.24000000033330823  \\
                0.25  0.25000000033330827  \\
                0.26  0.26000000033331516  \\
                0.27  0.27000000033331517  \\
                0.28  0.2800000003333152  \\
                0.29  0.29000000033331513  \\
                0.3  0.30000000033331514  \\
                0.31  0.31000000033331515  \\
                0.32  0.3200000003333083  \\
                0.33  0.3300000003333083  \\
                0.34  0.3400000003333049  \\
                0.35  0.3500000003332981  \\
                0.36  0.360000000333288  \\
                0.37  0.3700000003332713  \\
                0.38  0.38000000033323356  \\
                0.39  0.3900000003331075  \\
                0.4  0.39999987697718614  \\
                0.41  0.4000000033331391  \\
                0.42  0.40000000333323954  \\
                0.43  0.40000000333327307  \\
                0.44  0.4000000033332898  \\
                0.45  0.40000000333329977  \\
                0.46  0.4000000033333065  \\
                0.47  0.4000000033333113  \\
                0.48  0.400000003333315  \\
                0.49  0.40000000333331787  \\
                0.5  0.40000000333332003  \\
                0.51  0.4000000033333219  \\
                0.52  0.4000000033333234  \\
                0.53  0.40000000333332475  \\
                0.54  0.40000000333332586  \\
                0.55  0.40000000333332697  \\
                0.56  0.40000000333332775  \\
                0.57  0.40000000333332864  \\
                0.58  0.4000000033333293  \\
                0.59  0.40000000333332986  \\
                0.6  0.40000000333333047  \\
                0.61  0.4000000033333309  \\
                0.62  0.40000000333333147  \\
                0.63  0.4000000033333319  \\
                0.64  0.4000000033333323  \\
                0.65  0.4000000033333328  \\
                0.66  0.40000000333333313  \\
                0.67  0.4000000033333334  \\
                0.68  0.4000000033333338  \\
                0.69  0.40000000333333413  \\
                0.7  0.4000000033333343  \\
                0.71  0.40000000333333463  \\
                0.72  0.40000000333333496  \\
                0.73  0.4000000033333353  \\
                0.74  0.4000000033333356  \\
                0.75  0.4000000033333358  \\
                0.76  0.400000003333336  \\
                0.77  0.4000000033333363  \\
                0.78  0.40000000333333663  \\
                0.79  0.4000000033333369  \\
                0.8  0.40000000333333713  \\
                0.81  0.40000000333333735  \\
                0.82  0.4000000033333377  \\
                0.83  0.40000000333333796  \\
                0.84  0.40000000333333824  \\
                0.85  0.40000000333333857  \\
                0.86  0.40000000333333885  \\
                0.87  0.4000000033333392  \\
                0.88  0.40000000333333957  \\
                0.89  0.40000000333333996  \\
                0.9  0.4000000033333404  \\
                0.91  0.4000000033333409  \\
                0.92  0.4000000033333413  \\
                0.93  0.40000000333334196  \\
                0.94  0.40000000333334257  \\
                0.95  0.4000000033333433  \\
                0.96  0.4000000033333442  \\
                0.97  0.4000000033333452  \\
                0.98  0.4000000033333465  \\
                0.99  0.40000000333334795  \\
                1.0  0.4000000033333501  \\
                1.01  0.40000000333335284  \\
                1.02  0.4000000033333566  \\
                1.03  0.40000000333336255  \\
                1.04  0.4000000033333728  \\
                1.05  0.4000000033333955  \\
                1.06  0.40000000333348573  \\
                1.07  0.40139581968281124  \\
                1.08  0.40558894383825345  \\
                1.09  0.4097905673889501  \\
                1.1  0.4140006974428861  \\
                1.11  0.4182193410400716  \\
                1.12  0.42244650515228305  \\
                1.13  0.42668219668275986  \\
                1.14  0.430926422465909  \\
                1.15  0.4351791892670075  \\
                1.16  0.43944050378190924  \\
                1.17  0.4437103726367545  \\
                1.18  0.44798880238768224  \\
                1.19  0.4522757995205423  \\
                1.2  0.45657137045061463  \\
                1.21  0.46087552152232675  \\
                1.22  0.4651882590089769  \\
                1.23  0.4695095891124583  \\
                1.24  0.4738395179629861  \\
                1.25  0.47817805161882837  \\
                1.26  0.48252519606603833  \\
                1.27  0.48688095721818964  \\
                1.28  0.4912453409161156  \\
                1.29  0.49561835292764966  \\
                1.3  0.49999999894736974  \\
                1.31  0.504390284596345  \\
                1.32  0.5087892154218858  \\
                1.33  0.5131967968972964  \\
                1.34  0.5176130344216303  \\
                1.35  0.5220379333194493  \\
                1.36  0.5264714988405855  \\
                1.37  0.5309137361599049  \\
                1.38  0.5353646503770758  \\
                1.39  0.5398242465163396  \\
                1.4  0.5442925295262849  \\
                1.41  0.5487695042796236  \\
                1.42  0.5532551755729718  \\
                1.43  0.5577495481266329  \\
                1.44  0.5622526265843844  \\
                1.45  0.5667644155132671  \\
                1.46  0.5712849194033777  \\
                1.47  0.5758141426676663  \\
                1.48  0.5803520896417345  \\
                1.49  0.5848987645836391  \\
                1.5  0.589454171673697  \\
            }
            ;
        \addlegendentry {$x_{1}$}
        \addplot[color=rojito, thick, solid]
            table[row sep={\\}]
            {
                \\
                0.01  -3.3331515150050646e-10  \\
                0.02  -3.333082446951048e-10  \\
                0.03  -3.333082446951048e-10  \\
                0.04  -3.333082446951048e-10  \\
                0.05  -3.333082446951048e-10  \\
                0.06  -3.3331515150050646e-10  \\
                0.07  -3.3331515150050646e-10  \\
                0.08  -3.3331515150050646e-10  \\
                0.09  -3.3330824469510474e-10  \\
                0.1  -3.3330824469510474e-10  \\
                0.11  -3.3330824469510474e-10  \\
                0.12  -3.3330824469510474e-10  \\
                0.13  -3.3330824469510474e-10  \\
                0.14  -3.3330824469510474e-10  \\
                0.15  -3.3330824469510474e-10  \\
                0.16  -3.3330824469510474e-10  \\
                0.17  -3.3330824469510474e-10  \\
                0.18  -3.3330824469510474e-10  \\
                0.19  -3.3330824469510474e-10  \\
                0.2  -3.3330824469510474e-10  \\
                0.21  -3.3330824469510474e-10  \\
                0.22  -3.3330824469510474e-10  \\
                0.23  -3.3330824469510474e-10  \\
                0.24  -3.3330824469510474e-10  \\
                0.25  -3.3330824469510474e-10  \\
                0.26  -3.3331515150050646e-10  \\
                0.27  -3.3331515150050646e-10  \\
                0.28  -3.3331515150050646e-10  \\
                0.29  -3.3331515150050646e-10  \\
                0.3  -3.3331515150050646e-10  \\
                0.31  -3.3331515150050646e-10  \\
                0.32  -3.3330824469510474e-10  \\
                0.33  -3.3330824469510464e-10  \\
                0.34  -3.3330489476427665e-10  \\
                0.35  -3.3329812514966675e-10  \\
                0.36  -3.3328800029208e-10  \\
                0.37  -3.332713191177141e-10  \\
                0.38  -3.332335486393886e-10  \\
                0.39  -3.331075003301786e-10  \\
                0.4  1.230228138638062e-7  \\
                0.41  0.009999996666860828  \\
                0.42  0.019999996666760445  \\
                0.43  0.029999996666726897  \\
                0.44  0.03999999666671024  \\
                0.45  0.04999999666670022  \\
                0.46  0.059999996666693534  \\
                0.47  0.06999999666668864  \\
                0.48  0.07999999666668504  \\
                0.49  0.08999999666668214  \\
                0.5  0.09999999666667998  \\
                0.51  0.10999999666667812  \\
                0.52  0.11999999666667659  \\
                0.53  0.1299999966666752  \\
                0.54  0.13999999666667418  \\
                0.55  0.14999999666667307  \\
                0.56  0.15999999666667228  \\
                0.57  0.16999999666667132  \\
                0.58  0.17999999666667071  \\
                0.59  0.1899999966666701  \\
                0.6  0.19999999666666943  \\
                0.61  0.20999999666666902  \\
                0.62  0.2199999966666686  \\
                0.63  0.22999999666666815  \\
                0.64  0.23999999666666766  \\
                0.65  0.24999999666666728  \\
                0.66  0.2599999966666669  \\
                0.67  0.26999999666666663  \\
                0.68  0.2799999966666663  \\
                0.69  0.28999999666666576  \\
                0.7  0.29999999666666566  \\
                0.71  0.30999999666666533  \\
                0.72  0.31999999666666495  \\
                0.73  0.32999999666666463  \\
                0.74  0.3399999966666644  \\
                0.75  0.34999999666666426  \\
                0.76  0.35999999666666394  \\
                0.77  0.3699999966666638  \\
                0.78  0.3799999966666634  \\
                0.79  0.38999999666666313  \\
                0.8  0.3999999966666629  \\
                0.81  0.4099999966666627  \\
                0.82  0.41999999666666227  \\
                0.83  0.42999999666666205  \\
                0.84  0.43999999666666173  \\
                0.85  0.44999999666666146  \\
                0.86  0.4599999966666612  \\
                0.87  0.4699999966666608  \\
                0.88  0.4799999966666604  \\
                0.89  0.48999999666666005  \\
                0.9  0.4999999966666596  \\
                0.91  0.5099999966666592  \\
                0.92  0.5199999966666587  \\
                0.93  0.5299999966666582  \\
                0.94  0.5399999966666573  \\
                0.95  0.5499999966666567  \\
                0.96  0.5599999966666558  \\
                0.97  0.5699999966666548  \\
                0.98  0.5799999966666535  \\
                0.99  0.589999996666652  \\
                1.0  0.5999999966666498  \\
                1.01  0.6099999966666472  \\
                1.02  0.6199999966666434  \\
                1.03  0.6299999966666374  \\
                1.04  0.6399999966666272  \\
                1.05  0.6499999966666045  \\
                1.06  0.6599999966665143  \\
                1.07  0.6686041803171887  \\
                1.08  0.6744110561617467  \\
                1.09  0.6802094326110499  \\
                1.1  0.685999302557114  \\
                1.11  0.6917806589599285  \\
                1.12  0.6975534948477171  \\
                1.13  0.70331780331724  \\
                1.14  0.7090735775340908  \\
                1.15  0.7148208107329925  \\
                1.16  0.7205594962180907  \\
                1.17  0.7262896273632453  \\
                1.18  0.7320111976123178  \\
                1.19  0.7377242004794576  \\
                1.2  0.7434286295493854  \\
                1.21  0.7491244784776732  \\
                1.22  0.754811740991023  \\
                1.23  0.7604904108875417  \\
                1.24  0.7661604820370138  \\
                1.25  0.7718219483811716  \\
                1.26  0.7774748039339618  \\
                1.27  0.7831190427818104  \\
                1.28  0.7887546590838843  \\
                1.29  0.7943816470723503  \\
                1.3  0.8000000010526304  \\
                1.31  0.8056097154036551  \\
                1.32  0.8112107845781142  \\
                1.33  0.8168032031027038  \\
                1.34  0.8223869655783699  \\
                1.35  0.8279620666805507  \\
                1.36  0.8335285011594146  \\
                1.37  0.8390862638400953  \\
                1.38  0.8446353496229241  \\
                1.39  0.8501757534836603  \\
                1.4  0.855707470473715  \\
                1.41  0.8612304957203762  \\
                1.42  0.8667448244270282  \\
                1.43  0.8722504518733669  \\
                1.44  0.8777473734156156  \\
                1.45  0.883235584486733  \\
                1.46  0.8887150805966224  \\
                1.47  0.8941858573323337  \\
                1.48  0.8996479103582655  \\
                1.49  0.9051012354163609  \\
                1.5  0.910545828326303  \\
            }
            ;
        \addlegendentry {$x_{2}$}
        \addplot[color=black, dashed, forget plot]
            table[row sep={\\}]
            {
                \\
                0.01  0.01  \\
                1.5  1.5  \\
            };
    \end{axis}
    \end{tikzpicture}
    \caption{Rate splitting for the socially optimal (left) and selfish (right) routing policies as a function of input rate, for Example \ref{example:poa}.}
    \label{fig.poarates}
\end{figure}

Now consider the plot (on the right) of equilibrium rates for the selfish routing case. At low loads, there is no congestion and the solution is socially optimal. After queue 1 reaches capacity, selfish routing will continue to prefer station 1 until its queueing delay $\mu_1$ equates the difference $\kappa_2 -\kappa_1$ in transport delays. This happens at $r:=r_0\approx 0.4$; excess rates beyond this value will start choosing station 2 as in the socially optimal solution; but the waiting cost $\mu_1\cdot r_0$ implies inefficiency. Station 2 congests 
at $r=r_0+2/3$, after which selfish routing spreads the additional increase between both stations; this is required to maintain the indifference condition $\kappa_1+\mu_1 = \kappa_2+\mu_2$, and again implies some inefficiency.

\begin{figure}
    \centering
    \begin{tikzpicture}
        \begin{axis}[
            width=0.6\textwidth,
            height=0.4\textwidth,
            legend style={font=\footnotesize},
            axis x line=center,
            axis y line=center,
            xtick={0,0.33,1},
            xticklabels={$0$,$1/3$, $1$},
            ytick={0,10,20,30,40},
            xlabel={$r$},
            xlabel style={below},
            ylabel style={left},
            xmin=0,
            xmax=1.5,
            ymin=0,
            ymax = 45,
            ylabel={Cost},
            grid=both,
            legend pos = north west,
            legend cell align={left},]
            \addplot[color=azulcito, thick, solid]
                table[row sep={\\}]
                {
                    \\
                    0.01  0.010002575773304687  \\
                    0.02  0.020003129559159512  \\
                    0.03  0.03000373020733858  \\
                    0.04  0.0400040773803361  \\
                    0.05  0.05000433977606641  \\
                    0.06  0.060002431211560646  \\
                    0.07  0.07000244516113863  \\
                    0.08  0.08000243518782618  \\
                    0.09  0.09000443795309258  \\
                    0.1  0.10000437311103849  \\
                    0.11  0.11000437940138315  \\
                    0.12  0.12000436632418231  \\
                    0.13  0.13000434507846984  \\
                    0.14  0.14000432007446675  \\
                    0.15  0.15000245561238965  \\
                    0.16  0.1600024244523598  \\
                    0.17  0.17000442126895218  \\
                    0.18  0.18000438402131236  \\
                    0.19  0.190004396457117  \\
                    0.2  0.20000423470658177  \\
                    0.21  0.2100042142442732  \\
                    0.22  0.22000419611200853  \\
                    0.23  0.2300041806809045  \\
                    0.24  0.24000419474635978  \\
                    0.25  0.25000421026323844  \\
                    0.26  0.2600042253973103  \\
                    0.27  0.27000267894695384  \\
                    0.28  0.28000267768698567  \\
                    0.29  0.2900026684363694  \\
                    0.3  0.3000026262651268  \\
                    0.31  0.3100026717584051  \\
                    0.32  0.3200024486069646  \\
                    0.33  0.33000269636629465  \\
                    0.34  0.7399998070057243  \\
                    0.35  1.3499998070017771  \\
                    0.36  1.9599998070005307  \\
                    0.37  2.5699998069995384  \\
                    0.38  3.179999806997543  \\
                    0.39  3.7899998069753393  \\
                    0.4  4.000005933654101  \\
                    0.41  4.100003783066215  \\
                    0.42  4.200003695224517  \\
                    0.43  4.3000057577656765  \\
                    0.44  4.400005731007865  \\
                    0.45  4.500005687308018  \\
                    0.46  4.600005640048636  \\
                    0.47  4.700005602003024  \\
                    0.48  4.8000055760310065  \\
                    0.49  4.900005562305665  \\
                    0.5  5.00000554926553  \\
                    0.51  5.100005539666918  \\
                    0.52  5.200005532609536  \\
                    0.53  5.300005526583415  \\
                    0.54  5.40000552140234  \\
                    0.55  5.500005517117412  \\
                    0.56  5.600005518679401  \\
                    0.57  5.700005518545093  \\
                    0.58  5.800005517259901  \\
                    0.59  5.900005515113429  \\
                    0.6  6.000005512254419  \\
                    0.61  6.10000550875192  \\
                    0.62  6.2000055046307665  \\
                    0.63  6.300005499893807  \\
                    0.64  6.400005494536623  \\
                    0.65  6.500005488557784  \\
                    0.66  6.600005481966026  \\
                    0.67  6.700005474785403  \\
                    0.68  6.800005467058465  \\
                    0.69  6.90000545884806  \\
                    0.7  7.000005450237597  \\
                    0.71  7.100005441329921  \\
                    0.72  7.200005432244796  \\
                    0.73  7.300005423115001  \\
                    0.74  7.40000541419347  \\
                    0.75  7.500005405532235  \\
                    0.76  7.600005397121803  \\
                    0.77  7.700005389096214  \\
                    0.78  7.800005381569942  \\
                    0.79  7.900005374630075  \\
                    0.8  8.000005368329782  \\
                    0.81  8.100005362684024  \\
                    0.82  8.200005357668392  \\
                    0.83  8.300005353222312  \\
                    0.84  8.400005349257267  \\
                    0.85  8.500005345670699  \\
                    0.86  8.600005342365298  \\
                    0.87  8.700005066432416  \\
                    0.88  8.800005063365557  \\
                    0.89  8.90000506056678  \\
                    0.9  9.000005058186407  \\
                    0.91  9.100005056458087  \\
                    0.92  9.200005055669422  \\
                    0.93  9.30000505611008  \\
                    0.94  9.400005058009407  \\
                    0.95  9.500004954086759  \\
                    0.96  9.600004937095335  \\
                    0.97  9.700005783821045  \\
                    0.98  9.800005776511517  \\
                    0.99  9.9000057707502  \\
                    1.0  10.000005766832553  \\
                    1.01  10.100005765149946  \\
                    1.02  10.20000576616895  \\
                    1.03  10.30000577012555  \\
                    1.04  10.400005776537514  \\
                    1.05  10.500005787070483  \\
                    1.06  10.666256604578262  \\
                    1.07  11.329404733862345  \\
                    1.08  11.992484731776234  \\
                    1.09  12.655496540763624  \\
                    1.1  13.318440103778201  \\
                    1.11  13.981315364198892  \\
                    1.12  14.644122265830003  \\
                    1.13  15.306860752902894  \\
                    1.14  15.969530770077846  \\
                    1.15  16.632132262445825  \\
                    1.16  17.294665175530497  \\
                    1.17  17.95712945528995  \\
                    1.18  18.619525048118724  \\
                    1.19  19.28185190084948  \\
                    1.2  19.944109960755043  \\
                    1.21  20.606299175550006  \\
                    1.22  21.268419493392738  \\
                    1.23  21.93047086288708  \\
                    1.24  22.59245323308418  \\
                    1.25  23.254366553484232  \\
                    1.26  23.916210774038294  \\
                    1.27  24.577985845150025  \\
                    1.28  25.239691717677346  \\
                    1.29  25.901328342934328  \\
                    1.3  26.562895672692754  \\
                    1.31  27.224393659183896  \\
                    1.32  27.88582225510018  \\
                    1.33  28.54718141359688  \\
                    1.34  29.208471088293756  \\
                    1.35  29.869691233276683  \\
                    1.36  30.53084180309933  \\
                    1.37  31.191922752784762  \\
                    1.38  31.852934037826955  \\
                    1.39  32.513875614192585  \\
                    1.4  33.17474743832238  \\
                    1.41  33.835549467132815  \\
                    1.42  34.49628165801762  \\
                    1.43  35.156943968849305  \\
                    1.44  35.81753635798064  \\
                    1.45  36.47805878424623  \\
                    1.46  37.138511206964  \\
                    1.47  37.79889358593649  \\
                    1.48  38.459205881452576  \\
                    1.49  39.1194480542887  \\
                    1.5  39.77962006571036  \\
                }
                ;
            \addlegendentry {$C_s$}
            \addplot[color=rojito, thick, solid]
                table[row sep={\\}]
                {
                    \\
                    0.01  0.009999977000375205  \\
                    0.02  0.01999997700037474  \\
                    0.03  0.029999977000375513  \\
                    0.04  0.03999997700037441  \\
                    0.05  0.04999997700037489  \\
                    0.06  0.05999997700037463  \\
                    0.07  0.06999997700037437  \\
                    0.08  0.07999997700037408  \\
                    0.09  0.08999997700037551  \\
                    0.1  0.09999997700037504  \\
                    0.11  0.10999997700037442  \\
                    0.12  0.11999997700037449  \\
                    0.13  0.12999997700037452  \\
                    0.14  0.13999997700037445  \\
                    0.15  0.14999997700037446  \\
                    0.16  0.15999997700037558  \\
                    0.17  0.16999997700037522  \\
                    0.18  0.17999997700037487  \\
                    0.19  0.18999997700037463  \\
                    0.2  0.1999999770003744  \\
                    0.21  0.20999997700037418  \\
                    0.22  0.21999997700037427  \\
                    0.23  0.22999997700037417  \\
                    0.24  0.2399999770003741  \\
                    0.25  0.24999997700037405  \\
                    0.26  0.259999977000374  \\
                    0.27  0.26999997700037404  \\
                    0.28  0.2799999770003757  \\
                    0.29  0.289999977000374  \\
                    0.3  0.29999997700037345  \\
                    0.31  0.30999997700037285  \\
                    0.32  0.3199999770003731  \\
                    0.33  0.3299999770003736  \\
                    0.34  0.39999995150037687  \\
                    0.35  0.4999999515003762  \\
                    0.36  0.5999999515003764  \\
                    0.37  0.6999999515003765  \\
                    0.38  0.7999999515003767  \\
                    0.39  0.8999999515003764  \\
                    0.4  0.9999999515003762  \\
                    0.41  1.099999951500376  \\
                    0.42  1.1999999515003765  \\
                    0.43  1.2999999515003766  \\
                    0.44  1.3999999515003763  \\
                    0.45  1.4999999515003763  \\
                    0.46  1.5999999515003764  \\
                    0.47  1.699999951500376  \\
                    0.48  1.7999999515003766  \\
                    0.49  1.8999999515003763  \\
                    0.5  1.9999999515003763  \\
                    0.51  2.0999999515003758  \\
                    0.52  2.199999951500376  \\
                    0.53  2.2999999515003764  \\
                    0.54  2.399999951500377  \\
                    0.55  2.499999951500377  \\
                    0.56  2.599999951500377  \\
                    0.57  2.699999951500376  \\
                    0.58  2.799999951500376  \\
                    0.59  2.899999951500375  \\
                    0.6  2.9999999515003752  \\
                    0.61  3.0999999515003758  \\
                    0.62  3.199999951500376  \\
                    0.63  3.2999999515003764  \\
                    0.64  3.3999999515003774  \\
                    0.65  3.4999999515003766  \\
                    0.66  3.5999999515003767  \\
                    0.67  3.6999999515003767  \\
                    0.68  3.799999951500377  \\
                    0.69  3.899999951500376  \\
                    0.7  3.9999999515003752  \\
                    0.71  4.099999951500377  \\
                    0.72  4.199999951500375  \\
                    0.73  4.299999951500376  \\
                    0.74  4.399999951500376  \\
                    0.75  4.499999951500376  \\
                    0.76  4.599999951500377  \\
                    0.77  4.699999951500376  \\
                    0.78  4.799999951500376  \\
                    0.79  4.8999999515003765  \\
                    0.8  4.999999951500377  \\
                    0.81  5.099999951500377  \\
                    0.82  5.199999951500376  \\
                    0.83  5.299999951500376  \\
                    0.84  5.3999999515003765  \\
                    0.85  5.499999951500375  \\
                    0.86  5.599999951500376  \\
                    0.87  5.699999951500376  \\
                    0.88  5.799999951500377  \\
                    0.89  5.899999951500376  \\
                    0.9  5.999999951500376  \\
                    0.91  6.0999999515003775  \\
                    0.92  6.199999951500376  \\
                    0.93  6.299999951500378  \\
                    0.94  6.399999951500376  \\
                    0.95  6.499999951500377  \\
                    0.96  6.599999951500376  \\
                    0.97  6.699999951500375  \\
                    0.98  6.799999951500376  \\
                    0.99  6.899999951500376  \\
                    1.0  6.999999951500376  \\
                    1.01  7.609999458500381  \\
                    1.02  8.219999458500386  \\
                    1.03  8.82999945850038  \\
                    1.04  9.439999458500385  \\
                    1.05  10.04999945850038  \\
                    1.06  10.659999458500382  \\
                    1.07  11.269999458500385  \\
                    1.08  11.879999458500377  \\
                    1.09  12.489999458500387  \\
                    1.1  13.099999458500386  \\
                    1.11  13.709999458500391  \\
                    1.12  14.319999458500382  \\
                    1.13  14.929999458500381  \\
                    1.14  15.539999458500375  \\
                    1.15  16.149999458500375  \\
                    1.16  16.759999458500378  \\
                    1.17  17.369999458500374  \\
                    1.18  17.97999945850038  \\
                    1.19  18.58999945850038  \\
                    1.2  19.19999945850038  \\
                    1.21  19.80999945850038  \\
                    1.22  20.41999945850037  \\
                    1.23  21.029999458500384  \\
                    1.24  21.639999458500377  \\
                    1.25  22.249999458500376  \\
                    1.26  22.85999945850038  \\
                    1.27  23.46999945850039  \\
                    1.28  24.079999458500385  \\
                    1.29  24.68999945850038  \\
                    1.3  25.299999458500384  \\
                    1.31  25.90999945850038  \\
                    1.32  26.51999945850039  \\
                    1.33  27.129999458500386  \\
                    1.34  27.739999458500385  \\
                    1.35  28.34999945850039  \\
                    1.36  28.959999458500388  \\
                    1.37  29.569999458500384  \\
                    1.38  30.179999458500376  \\
                    1.39  30.789999458500375  \\
                    1.4  31.399999458500368  \\
                    1.41  32.009999458500374  \\
                    1.42  32.619999458500374  \\
                    1.43  33.22999945850038  \\
                    1.44  33.839999458500365  \\
                    1.45  34.44999945850038  \\
                    1.46  35.05999945850037  \\
                    1.47  35.669999458500385  \\
                    1.48  36.27999945850038  \\
                    1.49  36.889999458500384  \\
                    1.5  37.49999945850039  \\
                }
                ;
            \addlegendentry {$C^{opt}_s$}
            \addplot[color=verdecito, thick, solid]
                table[row sep={\\}]
                {
                    \\
                    0.01  0.009999997000348303  \\
                    0.02  0.019999997000496355  \\
                    0.03  0.02999999700059874  \\
                    0.04  0.03999999700066723  \\
                    0.05  0.04999999700072112  \\
                    0.06  0.059999997000330874  \\
                    0.07  0.06999999700033128  \\
                    0.08  0.07999999700032868  \\
                    0.09  0.0899999970007424  \\
                    0.1  0.09999999700073016  \\
                    0.11  0.10999999700073017  \\
                    0.12  0.11999999700072801  \\
                    0.13  0.12999999700072273  \\
                    0.14  0.13999999700071755  \\
                    0.15  0.14999999700033484  \\
                    0.16  0.15999999700032685  \\
                    0.17  0.16999999700074125  \\
                    0.18  0.1799999970007337  \\
                    0.19  0.18999999700073317  \\
                    0.2  0.19999999700070034  \\
                    0.21  0.20999999700069483  \\
                    0.22  0.21999999700069092  \\
                    0.23  0.22999999700068893  \\
                    0.24  0.2399999970006917  \\
                    0.25  0.2499999970006925  \\
                    0.26  0.25999999700069576  \\
                    0.27  0.26999999700035926  \\
                    0.28  0.2799999970003608  \\
                    0.29  0.2899999970003614  \\
                    0.3  0.29999999700035296  \\
                    0.31  0.3099999970003615  \\
                    0.32  0.3199999970003318  \\
                    0.33  0.3299999970003653  \\
                    0.34  0.3439474475473731  \\
                    0.35  0.37419670504030733  \\
                    0.36  0.42077916094448065  \\
                    0.37  0.48279967267768603  \\
                    0.38  0.5594347267670409  \\
                    0.39  0.6499249946533124  \\
                    0.4  0.7496213800450713  \\
                    0.41  0.8496213800447772  \\
                    0.42  0.9496213800447613  \\
                    0.43  1.0496213800450456  \\
                    0.44  1.149621380045042  \\
                    0.45  1.2496213800450386  \\
                    0.46  1.3496213800450327  \\
                    0.47  1.449621380045024  \\
                    0.48  1.5496213800450171  \\
                    0.49  1.649621380045017  \\
                    0.5  1.749621380045017  \\
                    0.51  1.8496213800450148  \\
                    0.52  1.9496213800450168  \\
                    0.53  2.0496213800450174  \\
                    0.54  2.1496213800450095  \\
                    0.55  2.249621380045017  \\
                    0.56  2.349621380045015  \\
                    0.57  2.4496213800450097  \\
                    0.58  2.5496213800450085  \\
                    0.59  2.6496213800450077  \\
                    0.6  2.7496213800450158  \\
                    0.61  2.849621380045012  \\
                    0.62  2.9496213800450106  \\
                    0.63  3.049621380045011  \\
                    0.64  3.1496213800450072  \\
                    0.65  3.24962138004501  \\
                    0.66  3.3496213800450065  \\
                    0.67  3.4496213800450075  \\
                    0.68  3.5496213800450014  \\
                    0.69  3.649621380044999  \\
                    0.7  3.7496213800450016  \\
                    0.71  3.8496213800450017  \\
                    0.72  3.9496213800450017  \\
                    0.73  4.0496213800449965  \\
                    0.74  4.149621380044997  \\
                    0.75  4.249621380045  \\
                    0.76  4.349621380044997  \\
                    0.77  4.449621380044993  \\
                    0.78  4.549621380044988  \\
                    0.79  4.649621380044991  \\
                    0.8  4.749621380044994  \\
                    0.81  4.849621380044995  \\
                    0.82  4.94962138004499  \\
                    0.83  5.049621380044987  \\
                    0.84  5.149621380044989  \\
                    0.85  5.249621380044987  \\
                    0.86  5.349621380044993  \\
                    0.87  5.44962138004495  \\
                    0.88  5.54962138004495  \\
                    0.89  5.649621380044953  \\
                    0.9  5.7496213800449505  \\
                    0.91  5.849621380044953  \\
                    0.92  5.949621380044953  \\
                    0.93  6.0496213800449565  \\
                    0.94  6.149621380044954  \\
                    0.95  6.24962138004494  \\
                    0.96  6.349621380044929  \\
                    0.97  6.449621380045054  \\
                    0.98  6.5496213800450525  \\
                    0.99  6.649621380045042  \\
                    1.0  6.749621380045046  \\
                    1.01  6.849621380045045  \\
                    1.02  6.949621380045043  \\
                    1.03  7.049621380045047  \\
                    1.04  7.149621380045047  \\
                    1.05  7.249621380045052  \\
                    1.06  7.349658161817469  \\
                    1.07  7.452920117492984  \\
                    1.08  7.561390075804959  \\
                    1.09  7.674971406428128  \\
                    1.1  7.793570125995153  \\
                    1.11  7.917094802279638  \\
                    1.12  8.045456462676068  \\
                    1.13  8.178568506747457  \\
                    1.14  8.31634662262574  \\
                    1.15  8.45870870706289  \\
                    1.16  8.605574788943255  \\
                    1.17  8.756866956078728  \\
                    1.18  8.912509285119135  \\
                    1.19  9.072427774420104  \\
                    1.2  9.236550279720039  \\
                    1.21  9.404806452486003  \\
                    1.22  9.577127680797112  \\
                    1.23  9.75344703264083  \\
                    1.24  9.933699201505078  \\
                    1.25  10.117820454155527  \\
                    1.26  10.305748580493486  \\
                    1.27  10.497422845395814  \\
                    1.28  10.692783942443548  \\
                    1.29  10.891773949451007  \\
                    1.3  11.094336285712014  \\
                    1.31  11.300415670884218  \\
                    1.32  11.50995808543685  \\
                    1.33  11.722910732591055  \\
                    1.34  11.939222001685753  \\
                    1.35  12.158841432905351  \\
                    1.36  12.381719683309134  \\
                    1.37  12.607808494105038  \\
                    1.38  12.837060659113446  \\
                    1.39  13.06942999436962  \\
                    1.4  13.304871308815605  \\
                    1.41  13.543340376035182  \\
                    1.42  13.78479390698768  \\
                    1.43  14.029189523698433  \\
                    1.44  14.276485733866078  \\
                    1.45  14.526641906348454  \\
                    1.46  14.779618247490989  \\
                    1.47  15.035375778262946  \\
                    1.48  15.293876312168866  \\
                    1.49  15.555082433903621  \\
                    1.5  15.818957478721487  \\
                }
                ;
            \addlegendentry {$C_0$}
        \end{axis}
        \end{tikzpicture}
    \caption{Cost comparison for the socially optimal and selfish routing policies.}
    \label{fig.poa}
\end{figure}
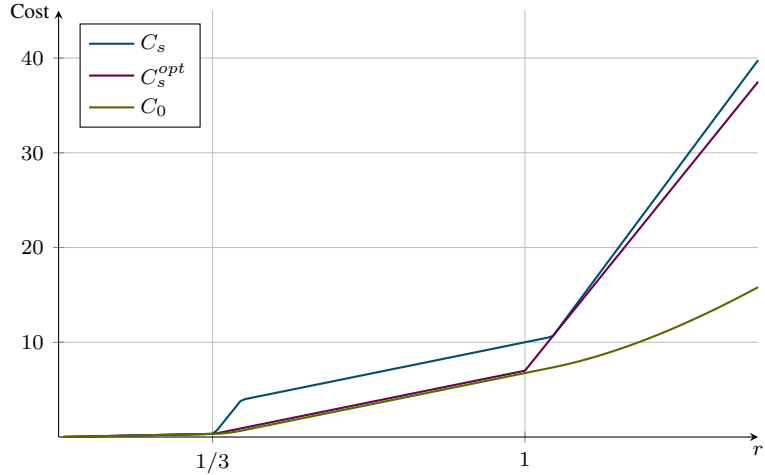

Fig. \ref{fig.poa} shows, as a function of $r$: the cost $C_0$ characterizing the selfish equilibrium, the resulting social cost $C_s$, and the optimum social cost $C_s^{opt}$. There is no price of anarchy at low loads; in the intermediate region inefficiency appears, for the most part constant at $(\kappa_2-\kappa_1)r_0$. There is a single point of efficiency again when both stations reach congestion, after which we observe a roughly linear gap in cost (constant price of anarchy) between both solutions. 
\end{example}
	
\subsection{Convergence}\label{ssec.conv}

Beyond the characterization of the equilibrium, we will establish that it is globally attractive under the dynamics \eqref{eq.sojourn dynamics}. From the Lipschitz nature of the field we know that from any initial condition $q(0)$, the trajectory $q(t)$ is well-defined for all time, and continuously differentiable. From this we can introduce the functions $\mu(q(t))$ and $D(\mu(q(t)))$; a key step of the convergence analysis will be to prove the latter is monotonically increasing along trajectories. 

These composite functions are not, however, differentiable. In particular the formula $\mu_j(q_j)=T \left[1-{c_j}/{q_j}\right]^+$ is not differentiable at $q_j=c_j$ (see Fig. \ref{fig.mu}). Similarly, if we compose it with $D_{2j}(\mu_j)$ in \eqref{eq.d2j} we obtain:
\begin{align}\label{eq.d2jqj}
D_{2j}(\mu_j(q_j)) 
= \begin{cases}
 0, &  q_j\leq c_j;\\
 c_j \log\left(\frac{c_j}{q_j}\right), & q_j > c_j;
\end{cases}
\end{align}
this (non-positive) function is not differentiable at $q_j=c_j$.

Nevertheless, both  $\mu_j(q_j)$ and $D_{2j}(\mu_j(q_j))$ are Lipschitz functions; composing them with the smooth $q(t)$ will yield \emph{absolutely continuous} functions $\mu_j(q_j(t))$ and $D_{2j}(\mu_j(q(t)))$; time derivatives exist almost everywhere, and the functions are integrals in time of their derivatives. We will denote by
$\dot{\mu}_j$ and $\dot{D}_{2j}$ these derivatives, and let 
$\mathcal{T}\subset \mathbb{R}_+$ be the set of times for which they are well-defined for every $j$. The complement of $\mathcal{T}$ has Lebesgue measure zero.  

\begin{lemma}\label{lem.dots}
    For every $t\in \mathcal{T}$, we have 
    \begin{align}\label{eq.chain}        
    \dot{D}_{2j}(t)=\frac{-c_j}{T-\mu_j}\dot{\mu}_j(t);
    \end{align}
    furthermore, if 
    $q_j(t)\leq c_j$, then $\dot{\mu}_j=\dot{D}_{2j}=0$.
    \end{lemma}
Proof is given in Appendix \ref{app.lemma}. We now get to the monotonicity result. 

\begin{proposition} \label{prop.dinc}
Consider a trajectory $q(t)$ of \eqref{eq.sojourn dynamics}. For any $t\in \mathcal{T}$ defined above, the dual function \eqref{eq.dual} satisfies:
\begin{align}\label{eq.ddot}
		\frac{d}{dt}{D}(\mu(q(t))) &
		=  \sum_{j:q_j(t) > c_j}T c_j \left[\frac{\dot{q}_j}{q_j}\right]^2 \geq 0. 
\end{align} 
Consequently, $D(\mu(q))$ is non-decreasing along trajectories $q(t)$ arising from the dynamics \eqref{eq.sojourn dynamics}.
 \end{proposition}
 
\begin{IEEEproof} 
We differentiate separately each term of the dual function; since $D_1(\mu)$ is differentiable we apply the chain-rule and the gradient in \eqref{eq.d1grad} to obtain:
\[
\dot{D_1} = \ip{\nabla D_1(\mu)}{\dot{\mu}} = \sum_i \ip{x^i(\mu)}{\dot{\mu}} = \ip{\dot{q}+q/T}{\dot{\mu}};
\]
here $\ip{\cdot}{\cdot}$ is the standard inner product in $\mathbb{R}^n$, and we have invoked \eqref{eq.sojourn dynamics-state}. Note that $\dot{\mu}$ is well-defined within $\mathcal{T}$. 

Now consider $D_2(\mu)$, and apply Lemma \ref{lem.dots} to obtain:
\begin{equation*}
    \dot{D_2} = \sum_j \dot{D}_{2j}=\sum_j \frac{-c_j}{T-\mu_j}\dot{\mu}_j =\sum_{j:q_j(t) > c_j} \frac{-q_j}{T} \dot{\mu}_j 
 = -\ip{{q}/{T}}{\dot{\mu}}.
\end{equation*}
For the third step we have used two facts: first, according to Lemma \ref{lem.dots} we may ignore in the sum terms where $q_j(t)\leq c_j$, since $\dot{\mu_j}=0$ in that case. Second, for $j$ where $q_j(t)>c_j$, we substitute $\mu_j=T(1 - c_j/q_j)$ from \eqref{eq.sojourn dynamics-mu}.

Adding both derivatives and cancelling terms we find
\begin{equation*}
\dot{D} = \ip{\dot{q}}{\dot{\mu}} = \sum_j \dot{q}_j \dot{\mu_j} =\sum_{j:q_j(t) > c_j} 
\dot{q}_j  \mu'_j(q_j) \dot{q}_j =\sum_{j:q_j(t) > c_j} \frac{Tc_j}{q_j^2} [\dot{q}_j]^2;\quad t \in \mathcal{T}.
\end{equation*}
Again we have removed terms with $q_j\leq c_j$, and in the rest we may apply the chain rule. This establishes \eqref{eq.ddot} almost everywhere.  Integrating in time and invoking absolute continuity, $D$ must be non-decreasing along trajectories. 
 \end{IEEEproof}

Since $D(\mu(q(t)))$ is non-decreasing and bounded above, it must approach a finite limit as $t\to \infty$. We will show that this limit must be the global maximum $D(\mu^*)$, and our trajectory  converges to the unique equilibrium. We state the corresponding result. 

\begin{theorem}\label{teo.convergence}
Any trajectory $(q(t),X(t),\mu(t))$ of \eqref{eq.sojourn dynamics} converges asymptotically to the equilibrium point $(X^*,q^*,\mu^*)$.    
\end{theorem}

The proof is based on a Lyapunov-LaSalle type argument \cite{khalil2002nonlinear} with the function $D(\mu(q))$, exploiting \eqref{eq.ddot} to characterize points where $\dot{D}=0$. Some non-standard modifications are required, details are given in Appendix \ref{app.theorem}.

\section{Elastic demand}\label{sec.elastic}

In this section we add a new ingredient to the formulation: instead of having a rigid quantity $r_i$ of EVs per second seeking charge from each location $i$, we consider the fact that demand itself is sensitive to the delays of the installation. 
Namely, if $\tau_i$ is the delay that EVs from location $i$ experience before receiving service, some customers may not be willing to wait this much and desist from  entering the system. 

To model this phenomenon, assume that arriving customers at location $i$ have a willingness-to-wait characteristic, which is a random variable $T_i$ with complementary cumulative distribution function $p_i(\tau_i) = \mathcal{P}(T_i > \tau_i)$. $T_i$ is assumed to be absolutely continuous, so $p_i(\tau_i)$ is continuous, decreasing from 1 to 0 as $\tau_i$ moves in the range $[0, \infty)$.

Randomness reflects customer heterogeneity; now, for a deterministic fluid model, we will not incorporate stochasticity at the individual EV level but rather consider its effect in large numbers. 
In this sense, if $\bar{r}_i$ is the (fixed) maximum rate of requests originating at  $i$, this rate will be ``thinned" by the probability $p_i(\tau_i)$, giving 
\begin{align}\label{eq.thinning}
    r_i(\tau_i) = \bar{r}_ip_i(\tau_i).
\end{align}
\begin{example}
 Assume that $T_i$ is a uniform random variable in the interval $[0,T]$, where $T$ is the maximum sojourn time. In that case we have the linear thinning 
    \[
    r_i(\tau_i) = \bar{r}_i\left[1 - \tau_i/T\right]^+,
    \]   
    shutting down demand as $\tau_i$ reaches $T$.
\end{example}

\subsection{Utility function representation}\label{ssec.util}

It will be useful to interpret \eqref{eq.thinning} as a \emph{demand curve} that specifies the quantity $r_i$ as a function of the ``price" $\tau_i$. In microeconomic language, this is equivalent to introducing an increasing, concave \emph{utility function} $U_i(r_i)$, and stating that the rate is chosen according to 
\begin{align}\label{eq.util argmax}
r_i = \arg\max_{r_i\geq 0} [U_i(r_i)-\tau_i r_i].    
\end{align}

Assuming differentiability, 
the above amounts to the first order condition
$U'_i(r_i) = \tau_i$; this should be the inverse function of \eqref{eq.thinning}. 
\begin{example}
Continuing with the uniform example, we find the corresponding utility.  $U'_i(r_i)$ must be linear, decreasing from $T$ to zero in the interval $[0,\bar{r}_i]$. Therefore:
    \[
    U_i(r_i)  = \begin{cases} T r_i\left(1 - \frac{r_i}{2\bar{r}_i}\right) & 0 \leq r_i \leq \bar{r}_i; \\
    T \frac{\bar{r}_i}{2}, & r_i > \bar{r}_i. 
    \end{cases}
    \]   
\end{example}

Note that utility will saturate (demand will satiate) at $r_i=\bar{r}_i$; to
simplify matters we will assume the following property, which holds in the example above: 
\begin{assumption}\label{ass.strict}
$U_i(r_i)$ is \emph{strictly} concave within the interval $r_i \in [0,\bar{r}_i]$.     
\end{assumption}

\subsection{Dynamics under elastic demand}\label{ssec.elastic dyn}

We now incorporate the elastic demand feature into the framework of Section \ref{sec.formu}.  We need to specify the appropriate $\tau_i$ for which to apply the expression \eqref{eq.thinning}. Under selfish routing to 
 the station(s) with smallest travel + waiting time, the delay-to-service experienced would be 
\[
\tau_i = \min_{j}(\kappa_{ij} + \mu_j) = \varphi(\kappa^i + \mu),
\]
invoking the notation \eqref{eq.varphi inicial}. Once again, for smoothness reasons we will use instead the log-sum-exp approximation \eqref{eq.logsumexp} to the minimum, namely:
\begin{align}\label{eq.tauimu}
\tau_i(\mu)= \varphi_\epsilon(\kappa^i + \mu) = \varphi_\epsilon^i(\mu),
\end{align}
recalling the notation \eqref{eq.logsumexp-i}. We are now ready to formulate the new dynamic model, incorporating into \eqref{eq.sojourn dynamics} the elastic rate thinning in \eqref{eq.thinning}:
\begin{subequations}\label{eq.elastic dynamics}
	\begin{align}
		\dot{q}_j &= \sum_{i=1}^m x_{ij} - \frac{q_j}{T}. \quad j=1,\ldots n.\label{eq.elastic dynamics-state}\\
		\mu_j(q_j) &= T \left[1-\frac{c_j}{q_j}\right]^+,  \quad j=1,\ldots n. \label{eq.elastic dynamics-mu} \\
		x_{ij} &=r_i \delta_{ij}(\mu), \ i=1,\ldots m, \ j=1,\ldots n,\mbox{ with } \delta_{ij}(\mu) \mbox{ in \eqref{eq.softargmin-i}};   \label{eq.elastic dynamics-x}\\
        r_i  &= \bar{r}_i p_i(\varphi_\epsilon^i(\mu)). \label{eq.elastic dynamics-r}
	\end{align}
\end{subequations}
As in the case of \eqref{eq.sojourn dynamics}, by substitution the above dynamics is equivalent to a differential equation $\dot{q}=f(q)$ with a globally Lipschitz right-hand side; the assumptions on $p_i(\cdot)$ preserve this feature. Hence, we again have unique solutions, defined for all time. We will also analyze the new dynamics through convex optimization.

\subsection{Equilibrium characterization} \label{ssec.elastic eq}

We denote here by $C(X,q,r)$ the cost function in \eqref{eq.costfc}; as compared to  Section \ref{sec.opt}, we now make explicit its dependence on the vector of demand rates $r$, no longer fixed. As domain for $r$ take the set $R: = \Pi_{i=1}^n [0,\bar{r}_i]$; zero rates are included, in this regard we refer to Remark \ref{rem.ri zero}.

$C(X,q,r)$ is \emph{not} a jointly convex function with this new variable; nevertheless, convexity holds for the solution of Problem \eqref{eq.inelastic opt} as a function of $r$. 
\begin{proposition}\label{prop.cr}
Let $\mathcal{C}(r)$ be the minimum of $C(X,q,r)$ over $(X,q)$ satisfying \eqref{eq.demand constr}-\eqref{eq.supply constr}
. Then $\mathcal{C}(r)$ is convex. 
\end{proposition}
\begin{IEEEproof}
    From our analysis of duality in the previous section, we know that 
    \[
    \mathcal{C}(r) = \max_\mu D(r,\mu),
    \]
    where $D(r,\mu)$ is the dual cost in \eqref{eq.dual}, making explicit its dependence on $r$. Noting the $D(r,\mu)$ is linear in $r$, the conclusion holds. 
\end{IEEEproof}

The optimal cost for a vector of demand rates is naturally  combined with the utility associated with these rates, defining the following convex objective:
\begin{align} \label{eq.deficit}
\psi(r): = \mathcal{C}(r) - \sum_i U_i(r_i).
\end{align}
Minimizing $\psi(r)$ is equivalent to maximizing $-\psi(r)$, a \emph{social suprlus} (utility minus cost); we will show below that our dynamics converges to precisely these optimal demand rates. 

To pursue the analysis, it is convenient to introduce 
\begin{equation}\label{eq.Wdef}
    W(r,\mu): = D(r,\mu) - \sum_i U_i(r_i),
 =\sum_i [r_i\varphi^{i}_\epsilon(\mu)-U_i(r_i)]+\sum_j  c_j \log\left(1-\frac{\mu_j}{T}\right), \nonumber 
\end{equation}
defined for $(r,\mu) \in R \times [0,T)^n$. Note the following properties:
\begin{itemize}
    \item For fixed $r$, $W(r,\mu)$ is strictly concave in $\mu\in [0,T)^n$. $W(r,\mu)\to -\infty$ when $\mu_j \uparrow T$. From Proposition \ref{prop.cr} we also obtain: 
    \[\psi(r) = \max_\mu  W(r,\mu).\]
    \item For fixed $\mu$, $W(r,\mu)$ is convex in $r\in R$; under Assumption \ref{ass.strict}, it is strictly convex. 
\end{itemize}
We now identify the saddle points of $W(r,\mu)$. 

\begin{proposition}\label{prop.saddle psi}
There exists a unique $(r^*,\mu^*) \in R\times [0,T)^n$, such that
\[
W(r^*,\mu) \leq W(r^*,\mu^*) \leq W(r,\mu^*) \ \ \forall  r\in R, \mu \in [0,T)^n.
\]
Furthermore, $r^* = \arg\min \psi(r)$.
\end{proposition}
\begin{IEEEproof}
Existence of a saddle point for a convex-concave function is a minimax result, of which there are multiple versions; in this case with bounded domains it can be found in \cite{rockafellar1970convex} (Corollary 37.6.1). 
Uniqueness of the saddle point follows from strict concavity/convexity. 
For the final statement, note:
\begin{equation*}
    \psi(r^*) =  \max_\mu  W(r^*,\mu) = W(r^*,\mu^*)  \leq W(r,\mu^*) \leq \max_\mu  W(r,\mu) = \psi(r). 
\end{equation*}
\end{IEEEproof}
We are now ready for our main result on equilibrium:


\begin{theorem}\label{teo.elastic-eq}
The following are equivalent: 
\begin{itemize}
\item[(i)] $(r^*,\mu^*)$ is the saddle point of  $W$ in \eqref{eq.Wdef}, and
given $r^*$,  $(X^*,q^*)$ are the solution to the optimization \eqref{eq.inelastic opt}. 
\item[(ii)] $(X^*,q^*,\mu^*,r^*)$ is an equilibrium point of 
\eqref{eq.elastic dynamics}.   
\end{itemize}
In particular, the dynamics have a unique equilibrium point.
\end{theorem}
\begin{IEEEproof}
Start from (i). At the saddle point $(r^*,\mu^*)$ we have $r^* = \arg\min_r W(r,\mu^*)$. From \eqref{eq.Wdef} we see that 
\begin{align*}
r^*_i & =  \arg\min_{r_i} [r_i\varphi^{i}_\epsilon(\mu^*)-U_i(r_i)] \\
 & = \arg\max_{r_i} [U_i(r_i)-r_i\varphi^{i}_\epsilon(\mu^*)] = \bar{r}_i p_i(\varphi^{i}_\epsilon(\mu^*)); 
\end{align*}
in the last step we used the characterization of the utility function in \eqref{eq.util argmax}. This equation is consistent with \eqref{eq.elastic dynamics-r}. Given $r^*_i$, we can invoke Theorem \ref{teo.sojourn-eq} to identify the saddle conditions for $(q^*,X^*,\mu^*)$ with the equilibrium of
\eqref{eq.elastic dynamics-state}-\eqref{eq.elastic dynamics-x}. Therefore we have (ii). 

Now start with (ii). For the given $r^*$, we know from Theorem \ref{teo.sojourn-eq} that $(q^*,X^*,\mu^*)$ are primal-dual optimal for the optimization \eqref{eq.inelastic opt}; in particular $\mu^*$ is dual optimal, 
\begin{align*}
    \mu^* &= \arg\max_\mu D(r^*,\mu) = 
\arg\max_\mu  W(r^*,\mu).
\end{align*}
Finally, \eqref{eq.elastic dynamics-r} implies $r^*=\arg\min_r  W(r,\mu^*)$ (these steps are reversible). Therefore $(r^*,\mu^*)$ is the saddle point of  $W$, and we have (i). 

Uniqueness of the equilibrium follows from the uniqueness of $(r^*,\mu^*)$ and the application of Theorem \ref{teo.sojourn-eq} for the remaining variables. 
\end{IEEEproof}

The previous result provides an interpretation of the equilibrium of the elastic dynamics, extending the one from Section \ref{ssec.eq-interpr}.
Indeed, if we consider jointly the two statements in condition (i), what the equilibrium achieves is the minimization of the cost $C(X,q,r) - \sum_i U_i(r_i)$ over all available degrees of freedom. Thus the dynamics achieve a specific welfare optimization, in which:
\begin{itemize}
    \item The utilities in question are directly related to the elastic thinning rule in \eqref{eq.thinning}.
    \item The cost portion $C(X,q,r)$ carries, again, the natural cost of transportation, a (negligible) regularization term, and the barrier terms $\beta_j(q_j)$. As in Section \ref{ssec.eq-interpr}, these exhibit a deviation with the respect to the waiting cost. Thus there may be inefficiency with respect to the natural social welfare optimization under elastic demand. 
    \end{itemize}

\subsection{Convergence}

Again, we show that our equilibrium is globally attractive under the dynamics 
\eqref{eq.elastic dynamics}. The analysis parallels Section \ref{ssec.conv}, based on the monotonicity of a certain function along trajectories of our dynamics. Here the function in question will be 
\begin{equation}
    \mathcal{D}(\mu) := \min_r W(r,\mu) = \mathcal{D}_1(\mu)+D_2(\mu) =  \sum_i  U^*_i(\varphi^{i}_\epsilon(\mu)) + \sum_j D_{2j}(\mu_j).\label{eq.elastic dual}  
\end{equation}
For the first term above we have introduced the Fenchel conjugate 
\begin{align}\label{eq.concave fenchel}    
U^*_i(\tau_i):=\min_{r_i}[\tau_i r_i-U_i(r_i)]
\end{align} 
of the concave utility function \cite{rockafellar1970convex}. The above minimum is achieved at $r_i(\tau_i) = \bar{r}_i p_i(\tau_i)$. Since the latter function is continuous, it is not hard to show that $U^*_i(\tau_i)$ is differentiable and non-decreasing in $\tau_i > 0$, with derivative $r_i(\tau_i)$. In \eqref{eq.elastic dual}, this function is composed with the smooth function $\varphi^{i}_\epsilon(\mu)$. Therefore the first term $\mathcal{D}_1$ in \eqref{eq.elastic dual} is differentiable in $\mu$, with 
\begin{align*} 
\nabla\mathcal{D}_1(\mu) = \sum_i \bar{r}_i p_i(\varphi^{i}_\epsilon(\mu)) 
\delta^i(\mu)
 = \sum_i x^i(\mu), 
    \end{align*}
where we have substituted expressions from the dynamics \eqref{eq.elastic dynamics}; this formula is completely analogous to \eqref{eq.d1grad} for the inelastic demand case. The second term $D_2(\mu)$  in \eqref{eq.elastic dual} is identical to \eqref{eq.dual}. Therefore we are in a position to replicate the analysis of
the previous section.

In particular, $\mathcal{D}(\mu(q(t)))$ will be absolutely continuous along trajectories of the dynamics, and we can identify a set $\mathcal{T}$ of times 
(whose complement has zero Lebesgue measure) where Lemma \ref{lem.dots} holds, as well as the following extension of Proposition \ref{prop.dinc}. The proof is analogous and is omitted. 

\begin{proposition} \label{prop.dinc elastic}
Consider a trajectory $q(t)$ of \eqref{eq.elastic dynamics}. For any $t\in \mathcal{T}$ defined above, the dual function \eqref{eq.elastic dual} satisfies:
\begin{align}\label{eq.ddote}
		\frac{d}{dt}{\mathcal{D}}(\mu(q(t))) &
		=  \sum_{j:q_j(t) > c_j}T c_j \left[\frac{\dot{q}_j}{q_j}\right]^2 \geq 0. 
\end{align} 
Consequently, ${\mathcal{D}}(\mu(q))$ is non-decreasing along trajectories $q(t)$ arising from the dynamics \eqref{eq.elastic dynamics}.
 \end{proposition}

Note, from the strict concavity of $W(r,\mu)$ in $\mu$, that $\mathcal{D}(\mu)$ in \eqref{eq.elastic dual} is strictly concave. Also, from the min-max inequality 
\begin{equation*}
 \max_\mu \mathcal{D}(\mu) = \max_\mu \min_r W(r,\mu) \leq \min_r\max_\mu W(r,\mu)= \min_r \psi(r) = \psi(r^*),
\end{equation*}
we conclude that $\mathcal{D}(\mu)$ is upper bounded by the optimal cost. There is equality above at the (unique) saddle point $(r^*,\mu^*)$, so we conclude that   $\mathcal{D}(\mu)$ achieves its maximum (only) at $\mu^*$. These are analogous properties to the ones we had for $D(\mu)$ in the previous section, leading to following convergence result. 

\begin{theorem}\label{teo.elastic convergence}
Any trajectory $(q(t),X(t),\mu(t),r(t))$ of \eqref{eq.elastic dynamics} converges asymptotically to the equilibrium point $(X^*,q^*,\mu^*,r^*)$ characterized in Theorem \ref{teo.elastic-eq}.    
\end{theorem}

The proof is very similar to the one for Theorem \ref{teo.convergence}, and is omitted. 
Some further remarks are  given in Appendix \ref{app.theorem}.
\section{Stochastic Simulations}\label{sec.sims}

In this section we present simulations of the selfish load balancing dynamics. A first objective is illustrating in a concrete geometric example the properties of the equilibrium (optimum of Problem \eqref{eq.inelastic opt}), and the convergence of trajectories.  A separate purpose is validating our fluid model as representative of more realistic conditions: in particular, we will simulate a \emph{stochastic} system generated by EVs arriving randomly in time and space, with random sojourn times, and stations assigned through selfish routing. We focus on the inelastic (fixed demand) case.

Our spatial domain is a square region, where recharge requests arrive as a Poisson process of overall rate $r=3$ EVs/min, spawning at a random spatial location with uniform distribution. Sojourn times are independent and exponentially distributed with mean $T=90$ min., so the total number of EVs present in steady-state is $rT=270$ on average. We fix $5$ charging stations
with $c_j=50$ slots each, at random positions. Since $rT>\sum_j c_j$, there is overall congestion; it will affect stations asymmetrically due to their locations. 

\begin{figure}[h]
	\centering
	\includegraphics[width=0.75\textwidth]{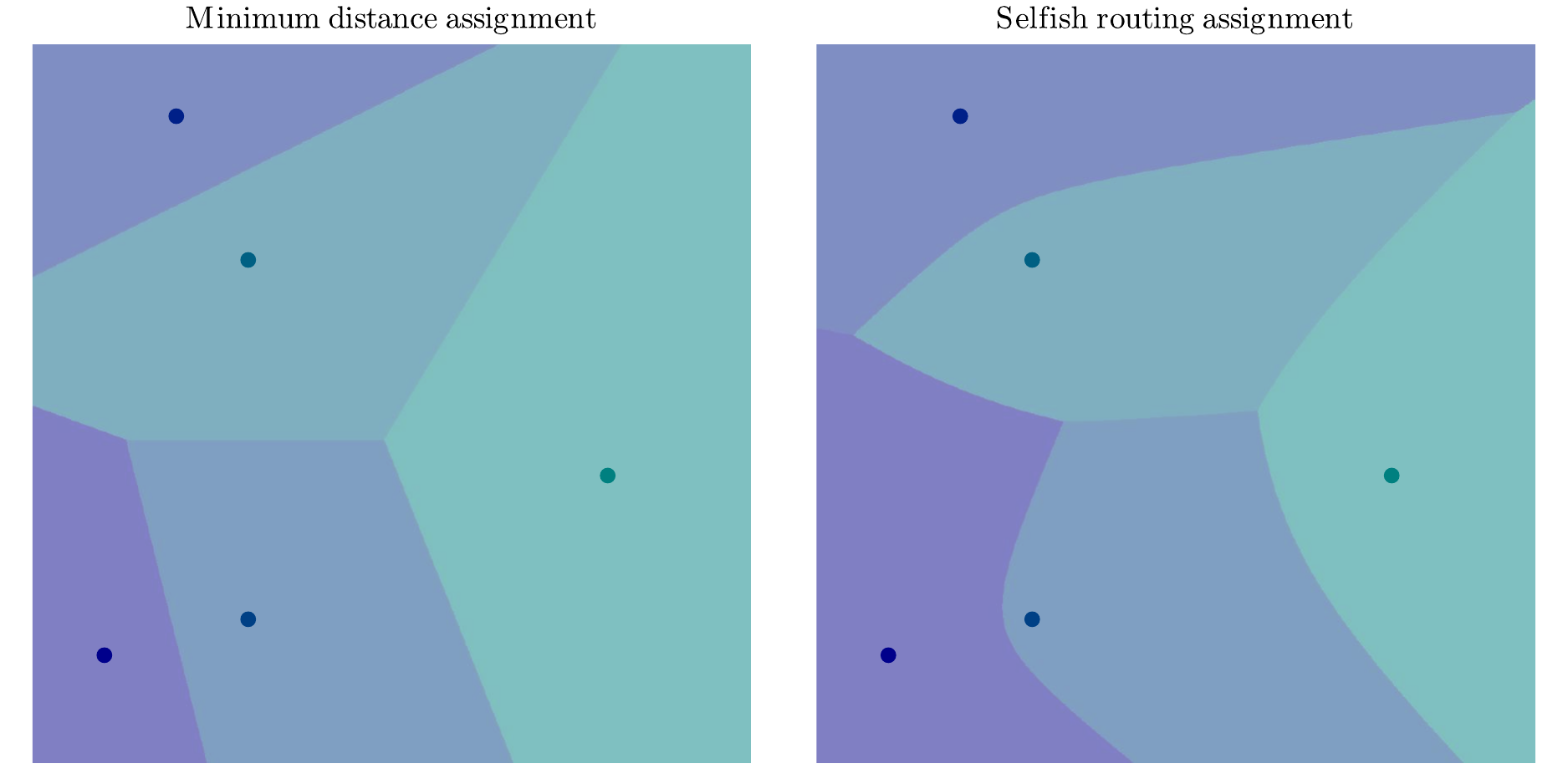}
	\caption{Charging station positions, minimum distance cells and attraction regions in equilibrium.} \label{fig.attraction}
\end{figure}

Travel times $\kappa_{ij}$ are modeled as the Euclidean distance divided by a speed $v$, chosen such that the travel time horizontally across the region is $50$ minutes. With the given station locations, the maximum time to reach the closest station is approximately $30$ min. On the left side of Fig. \ref{fig.attraction} we plot the station positions as well as the Voronoï cells that break up the plane if we use the closest station criterion. 

\begin{figure}
	\centering
	\includegraphics[width=0.6\textwidth]{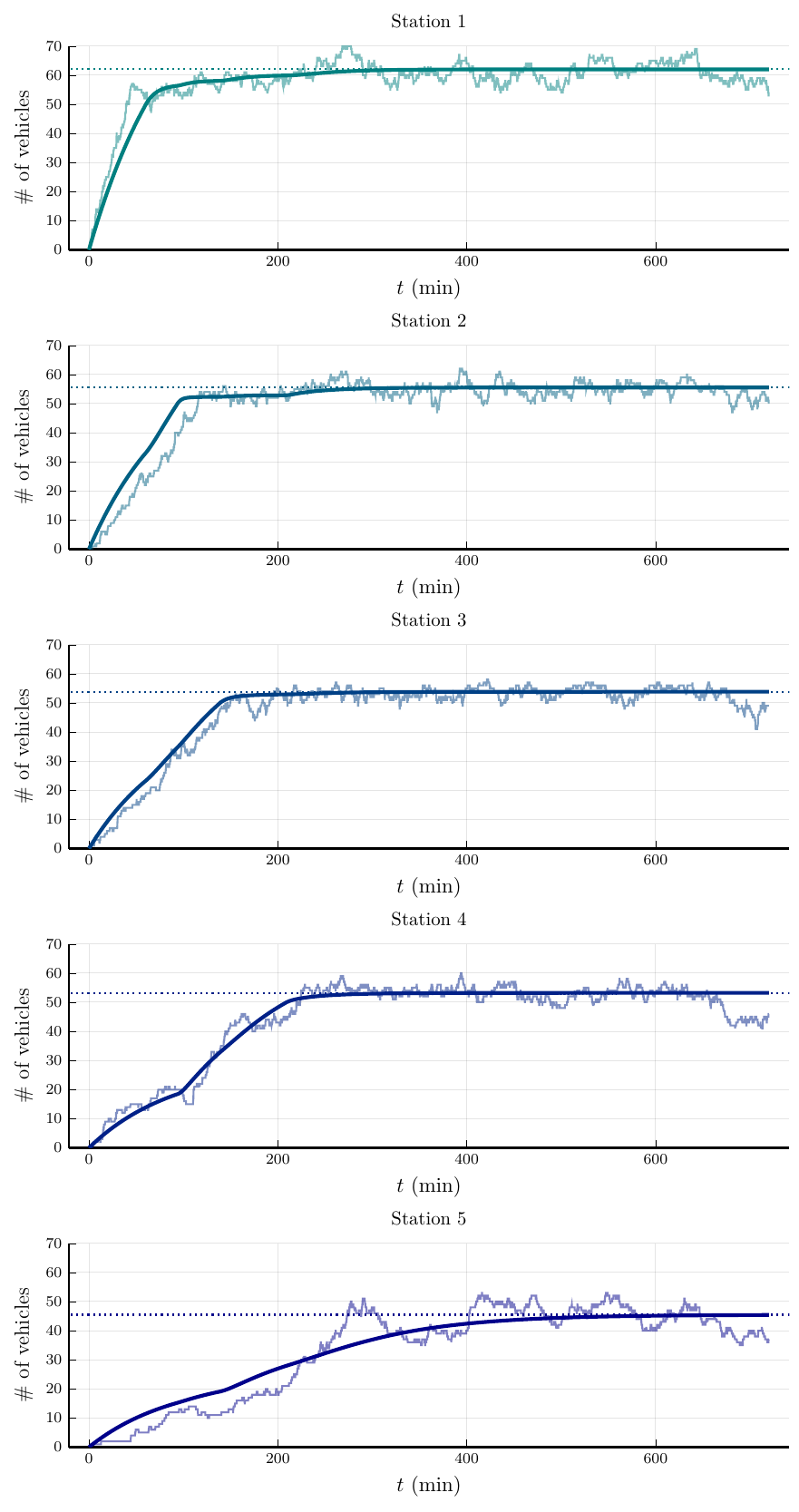}
	\caption{Time evolution of station occupations for the stochastic system and fluid approximation (thick lines). The predicted equilibrium of Theorem \ref{teo.sojourn-eq} is shown in dotted lines.} \label{fig.sojourn_time_sim}
\end{figure}

In Fig. \ref{fig.sojourn_time_sim} we plot the station occupation resulting from our spatial stochastic simulation, using the Julia library EVQueues \cite{evqueues.jl}. The selfish policy is applied, so that vehicles upon arrival are directed to the station that provides minimum delay to service. We also plot for comparison the computed numerical solutions of the fluid model \eqref{eq.sojourn dynamics}: since the model uses a finite set of arrival locations, we discretize the space to a uniform grid of $10000$ points. Note that the fluid model correctly captures both the transient and steady-state behavior of the stochastic system, and converges to the predicted equilibrium, the solution of Problem \eqref{eq.inelastic opt}, shown with dotted lines.

Initially, all stations are uncongested and EVs route themselves to the closest station. Station load is thus proportional to the area of the respective Voronoï cells, so it is asymmetric. In particular the easternmost station (labeled by 1) gets congested first: as its queueing delay builds up, neighboring stations start to receive more traffic; inflection points in the respective fluid trajectories signal this event. As the simulation progresses other stations reach congestion in succession, except for Station 5 (southwest point) which stays below capacity. 

In steady state, Station 1 operates with $q_1\approx 62$, i.e. $12$ requests in waiting, a queueing delay of $\mu_1\approx 17.4$ min. The remaining congested stations reach queueing delays $\mu_2\approx 8.9$, $\mu_3\approx 6.4$, and $\mu_4\approx 5.4$ minutes, while $\mu_5$ remains at $0$. 

The effect of these delays in steady-state routing is observed by plotting the attraction regions, on the right of Fig. \ref{fig.attraction}. 
Congested stations see their regions shrink, while less congested stations cover additional ground to compensate. Indifference curves can be shown to be arcs of hyperbolas in this case.

\section{Conclusions}\label{sec.concl}

 Network resource allocation calls for an interplay between dynamics and optimization. Besides helping design \emph{control} mechanisms for engineered systems with centrally adjudicated resources, optimization has also been applied successfully to characterize congestion \emph{games} between selfish agents. 

In this paper we have analyzed the operation of spatially distributed EV charging resources. Rather than a static, centrally planned optimal transport we have considered dynamic, selfish assignment by EV drivers endowed with delay information. We have proposed a dynamic model and analyzed its equilibrium and dynamics with tools of convex optimization. 

Natural lines for future research are: (i) Quantitative assessment of the price of anarchy, and its possible mitigation through more active control; 
(ii) Dynamics beyond constant demand (transients, tracking of daily variations); (iii) Stochastic versions of the queueing model and/or the geometry, for a more thorough analysis of variability.

\appendices

\section{Proof of Lemma \ref{lem.dots}}\label{app.lemma}
We must establish the identity \eqref{eq.chain} 
for functions $\mu_j(q_j(t))$ and $D_{2j}(\mu_j(q_j(t)))$ at times  when they are both differentiable. This includes all $t: q_j(t)\neq c_j$, since at these points
 $\mu_j(q_j)$ from \eqref{eq.mu of q}
 and $D_{2j}(\mu_j(q_j))$ from \eqref{eq.d2jqj} are both differentiable, and so is
 $q_j(t)$. In this case, \eqref{eq.chain} just follows from the chain rule,
 recalling that $D_{2j}(\mu_j) = c_j \log(1-\mu_j/T)$. 
Furthermore, if $q_j(t)<c_j$ both sides of \eqref{eq.chain} are zero. 

It remains to consider the case $q_j(t)=c_j$, where  $\mu_j(q_j)$ and $D_{2j}(\mu_j(q_j))$ are not smooth, they have finite but unequal lateral derivatives. Still, we are assuming $t\in \mathcal{T}$ so the composite functions $\mu_j(q_j(t))$ and $D_{2j}(\mu_j(q_j(t)))$ are differentiable. For this to happen their derivatives must be zero, as follows from a basic fact from differentiation, stated as Lemma  \ref{lem.derivadas laterales} below. So \eqref{eq.chain} holds in this case too. 
\hfill \rule{1.2ex}{1.2ex}

\begin{lemma}\label{lem.derivadas laterales}
Let $f:\mathbb{R}\to\mathbb{R}$ be continuous, $f(c)=0$ and with different lateral derivatives $f'(c-)\neq f'(c+)$. Let $q:\mathbb{R}\to\mathbb{R}$ be differentiable, $q(t_0)=c$. If the composite function  $g(t)=f(q(t))$ is differentiable at $t_0$, then $\dot{q}(t_0)=\dot{g}(t_0)=0$. 
\end{lemma}

\begin{IEEEproof}
Without loss of generality, take $t_0=0$, $c=0$. 
Write the first order Taylor expansion $q(t)=\alpha t + o(t)$.  

If $\alpha > 0$, then $q(t)$ has the same sign as $t$ in a neighborhood of zero, and thus:
\[
\lim_{t\to 0+} \frac{g(t)}{t} = 
\lim_{t\to 0+} \frac{f(q(t))}{q(t)}\frac{q(t)}{t}
 = 
f'(0+)\ \alpha.
\]
Analogously, $\lim_{t\to 0-} \frac{g(t)}{t} = f'(0-)\ \alpha$, a contradiction since $\dot{g}(0)$ exists and $f'(0-)\neq f'(0+)$. Similarly, we can rule out $\alpha < 0$. 
Therefore, $\alpha=\dot{q}(0)=0$. Now $q(t) = o(t)$ so 
\[
\lim_{t\to 0} \left| \frac{g(t)}{t}  \right|=\lim_{t\to 0}  \left|\frac{f(q(t))}{q(t)} \right|\left|\frac{q(t)}{t} \right| = 0;
\]
the first term is bounded because $f(q)$ has lateral derivatives. 
\end{IEEEproof}

\section{Proof of Theorem \ref{teo.convergence}} \label{app.theorem}

The function $V(q) = D(\mu(q))$ is well-defined over $q \in \mathbb{R}^n$, and bounded, with a maximum $V^* = D(\mu^*)$ achieved at $q=q^*$, the equilibrium queues. Furthermore, $V(q)$ is non-decreasing along trajectories $q(t)$ of the dynamics \eqref{eq.sojourn dynamics}. These elements suggest using $V$ as a Lyapunov function to establish convergence. 
A few difficulties arise, however: when some $\mu^*_j=0$, it is not true that $V(q)=V^*$ \emph{only} at equilibrium $q$. Furthermore, we have to deal with non-differentiability. Thus, we develop a specialized variant of the LaSalle principle \cite{khalil2002nonlinear}.

We begin by noting that the dynamics in $q(t)$ have a compact positively invariant set $[0,rT]^n$, where we denote $r=\sum_i r_i$. Indeed, since \eqref{eq.sojourn dynamics-x} implies $x_{ij}(t)< r_i \ \forall \ t$, we obtain from \eqref{eq.sojourn dynamics-state} the inequality 
\[
T\dot{q}_j < \sum_i Tr_i - q_{j} = T r- q_j.
\]
For $q_j \geq  rT$ we have $\dot{q}_j<0$; this implies a trajectory starting from $q_j\in [0,rT]$ cannot exceed the upper limit. In fact, a slight refinement of this argument implies that from any initial condition {outside} $[0,rT]^n$, this set is reached in finite time. Hence, it suffices to analyze the dynamics with initial conditions within this invariant set. 

Consider $q(0) \in [0,rT]^n$, and the resulting trajectory $q(t)$. From Proposition \ref{prop.dinc} we know that  $V(q(t))$ is monotonically non-decreasing, let its limit be $\bar{V}$. 

Let $L^+$ be the $\Omega$-limit set of the trajectory, itself an invariant set under the dynamics. By continuity, $V(q)\equiv \bar{V}$ for $q\in L^+$. 

Now consider an auxiliary trajectory, $\tilde{q}(t)$, with initial condition $q^+ \in L^+$. We conclude that 
\[
V(\tilde{q}(t))=D(\mu(\tilde{q}(t)))\equiv \bar{V} \ \Longrightarrow \ \frac{d}{dt}D (\mu(\tilde{q}(t)))\equiv 0.
\]
In reference to \eqref{eq.ddot}, this implies that $\dot{\tilde{q}}_j =0$ almost everywhere, for any $j: \tilde{q}_j(t)>c_j$; these congested queues must be at equilibrium, and will satisfy $\dot{\mu_j}\equiv 0$. But non-congested queues ($\tilde{q}_j(t)\leq c_j$) also satisfy  $\dot{\mu_j}= 0$ by Lemma \ref{lem.dots}; therefore $\mu(\tilde{q}(t))\equiv \tilde{\mu}$, constant. Consequently, the rates $X(\tilde{\mu})$ defined by \eqref{eq.sojourn dynamics-x} are also constant. 

Stations with $\tilde{\mu}_j = 0$ ($\tilde{q}_j(t)\leq c_j$) might not be at equilibrium; however they receive a constant input rate 
$\sum_i x_{ij}(\tilde{\mu})$ and evolve according to the first-order linear dynamics
\begin{align*}
	\dot{\tilde{q}}_j = \sum_i x_{ij}(\tilde{\mu}) - \tilde{q}_j/T. 
\end{align*}
This ODE has solutions that converge exponentially to $q_j^* = T \sum_i x_{ij}(\tilde{\mu})\in [0,c_j]$, an equilibrium. Thus, all queues reach equilibrium in the trajectory $\tilde{q}(t)$; since the equilibrium is unique from Theorem \ref{teo.sojourn-eq}, we conclude from this exercise that necessarily $\tilde{\mu} = \mu^*$, the dual optimal price, and $\bar{V} = V^* = D(\mu^*)$.

Return now to the \emph{original} trajectory $q(t)$. We know that
$V(q(t))=D(\mu(q(t))) \to D(\mu^*)$; since $D(\mu)$ is strictly concave, with a unique maximum, we must have $\mu(q(t)) \to \mu^*$.

Consequently, the input rate $\sum_{i}x_{ij}(\mu(q(t))$ to the $j$-th station converges to the optimal rate $\sum_{i}x^*_{ij}$. For the stable first order system \eqref{eq.sojourn dynamics-state}, it follows that the output 
$q_j(t)$ converges to the resulting equilibrium: $q_j(t) \to q^*_j$ for all $j$. 
\hfill \rule{1.2ex}{1.2ex}

\subsection*{Remark on the Proof of Theorem \ref{teo.elastic convergence}}
The preceding argument translates completely, replacing $D(\mu)$ in \eqref{eq.dual} by $\mathcal{D}(\mu)$ in \eqref{eq.elastic dual}. In the LaSalle-type argument one also obtains that $\mu(t)\equiv \tilde{\mu}$ constant, and 
it follows that the queue input rates are constant due to   \eqref{eq.elastic dynamics-x} and \eqref{eq.elastic dynamics-r}. Similarly, in the final step one obtains
 $\mu(t)\to \mu^*$ which implies input rates converging to a constant.


\end{document}